\documentclass[a4paper,11pt]{amsart} 
\usepackage{amscd,amssymb}
\usepackage{amsfonts}
\usepackage{graphicx,epsf}
\usepackage{enumerate,color}
\usepackage[mathcal]{euscript}
\usepackage{afterpage}
\numberwithin{equation}{section}

\usepackage{epstopdf}
\newcommand{\newtext}{}
\newcommand{\newmath}{}

\newcommand{\R}{\mathbb{R}}
\newcommand{\C}{\mathbb{C}}

\newcommand{\diver}{\mathrm{div}}

\DeclareMathOperator{\dbar}{\overline{\partial}}

\DeclareMathOperator{\bndry}{\partial\Omega}
\DeclareMathOperator{\T}{\mathbf{t}}

\newcommand{\argmin}{\mathop{\mathrm{arg\,min}}}
\newcommand{\sigmaCE}{\sigma^{\mbox{\tiny CE}}}
\newcommand{\sigmaAT}{\sigma^{\mbox{\tiny AT}}}

\begin{document}
\title[An Edge-Preserving D-bar Method for EIT]{A Data-Driven Edge-Preserving D-bar Method for Electrical Impedance Tomography}

\author[S. J. Hamilton, A. Hauptmann and S. Siltanen]{}


 \email{sarah.hamilton@marquette.edu}
 \email{andreas.hauptmann@helsinki.fi}
 \email{samuli.siltanen@helsinki.fi}


\begin{abstract}
In Electrical Impedance Tomography (EIT), the internal conductivity of a body is recovered via current and voltage measurements taken at its surface.  The reconstruction task is a highly ill-posed nonlinear inverse problem, which is very sensitive to noise, and requires the use of regularized solution methods, of which D-bar is the only proven method.  The resulting EIT images have low spatial resolution due to smoothing caused by low-pass filtered regularization.  In many applications, such as medical imaging, it is known \emph{a priori} that the target contains sharp features such as organ boundaries, as well as approximate ranges for realistic conductivity values.  In this paper,  we use this information in a new edge-preserving EIT algorithm, based on the original D-bar method coupled with a deblurring flow stopped at a minimal data discrepancy. The method makes heavy use of a novel data fidelity term based on the so-called {\em CGO sinogram}. This nonlinear data step provides superior robustness over traditional EIT data formats such as current-to-voltage matrices or Dirichlet-to-Neumann operators, \newtext{for commonly used current patterns}.
\end{abstract}

\maketitle

\centerline{\scshape Sarah Jane Hamilton }
\medskip
{\footnotesize
 \centerline{Department of Mathematics, Statistics, and Computer Science}
   \centerline{Marquette University}
   \centerline{ Milwaukee, WI 53233, USA}
} 

\medskip

\centerline{\scshape Andreas Hauptmann and Samuli Siltanen}
\medskip
{\footnotesize
 \centerline{ Department of Mathematics and Statistics }
   \centerline{University of Helsinki}
   \centerline{Helsinki, 00014, Finland}
}

\bigskip


\section{Introduction}
\begin{figure}[t]
\centering
\begin{picture}(325,110)
\put(0,0){\includegraphics[width=325 pt]{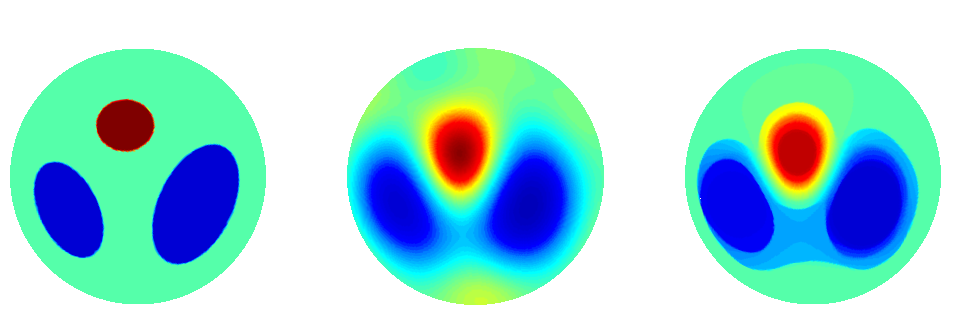}}
\put(30,102){\textbf{True}}
\put(142,102){\textbf{D-bar}}
\put(250,102){\textbf{Improved}}
\end{picture}
\caption{\label{fig:hammer_to_the_head}Left: true conductivity. Middle: blurry D-bar reconstruction from \newtext{0.5\% noise corrupted} EIT measurements, relative $l^1$-error 15.11\%. Right: reconstruction by the proposed edge-preserving algorithm, relative $l^1$-error 12.57\%.}
\end{figure}
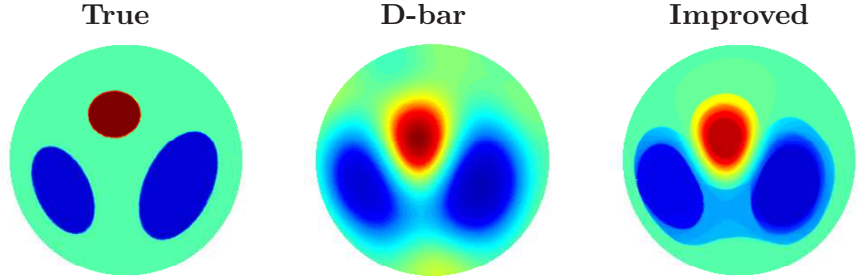

\noindent
Noise-robust solutions of \emph{ill-posed} inverse problems are based on regularization strategies.  For Electrical Impedance Tomography (EIT), the only proven regularization strategy is the \emph{low-pass filtered D-bar Method}, which sets high scattering frequencies to zero therefore resulting in smoothed images where sharp features important to applications such as medical imaging are often absent.  In this paper, we propose to recover edges in the smoothed EIT reconstructions by applying a deblurring flow stopped at a minimal data discrepancy (Figure~\ref{fig:hammer_to_the_head}), guided by a novel data fidelity term based on the so-called {\em CGO sinogram}, which provides superior robustness.

Electrical Impedance Tomography (EIT) is an imaging modality where an unknown physical body is probed with electricity using electrodes attached to the surface of the body. The goal is to recover the internal conductivity distribution of the body typically based on current-to-voltage boundary measurements. EIT has applications in medical imaging, underground contaminant detection and industrial process monitoring. See \cite{Cheney1999} and \cite[Chapter 12]{Mueller2012} for more details and applications of EIT.

We formulate the inverse conductivity problem \cite{Calder'on1980} for two-dimensional EIT in terms of voltage-to-current measurements. Let $\Omega\subset\R^2$ be the unit disc. We model the conductivity by a bounded measurable function $\sigma:\Omega\rightarrow\R$ satisfying $C>\sigma(z)\geq c>0$ for almost every $z\in\Omega$.  For a prescribed boundary voltage $f\in H^{1/2}(\bndry)$, the voltage potential $u$ satisfies the \emph{conductivity equation}
\begin{equation}\label{eq:cond_eq}
\begin{array}{rl}
\nabla\cdot\sigma\nabla u = & 0, \quad \mbox{ in }\Omega,\\
u|_{\bndry} =& \phi, \quad  \mbox{ on }\bndry.
\end{array}
\end{equation} 
Infinite-precision voltage and current measurements are modeled by the \emph{Dirichlet-to-Neumann} (DN), or voltage-to-current density, map
\begin{equation*}
\Lambda_\sigma:\;\phi\mapsto \left.\sigma\frac{\partial u}{\partial\nu}\right|_{\bndry},
\end{equation*}
where $\nu$ denotes the outward facing unit normal to $\bndry$. 
The goal of EIT is to recover the conductivity distribution $\sigma(z)$ for $z\in\Omega$, approximately, from the knowledge of a practical data matrix $\Lambda_\sigma^\delta$ satisfying $\|\Lambda_\sigma-\Lambda_\sigma^\delta\|_Y\leq \delta$ with known noise amplitude $\delta$ and an appropriate data space $Y$, see \cite{Knudsen2009} on details for such a space.

EIT is a severely ill-posed inverse problem, in fact it is only $\log$ stable.  By this we mean that small changes in boundary measurements can correspond to large changes in the internal conductivity distribution, and furthermore that noise in the data is amplified exponentially.  Therefore, {\em regularization} is needed for the noise-robust recovery of $\sigma$ from  $\Lambda_\sigma^\delta$. The forward map $\sigma\mapsto\Lambda_\sigma$ is too nonlinear to be covered by the presently available theory of iterative (Tikhonov-type) regularization. So far the only methodology providing proven regularization properties is the so-called {\em D-bar method} in dimension two \cite{Siltanen2000, Mueller2003,Knudsen2009}. Regularization for the 3D case is in progress, based on \cite{Cornean2006,Bikowski2010a,Delbary2011}.\footnote{Personal communication with Kim Knudsen and his team.}

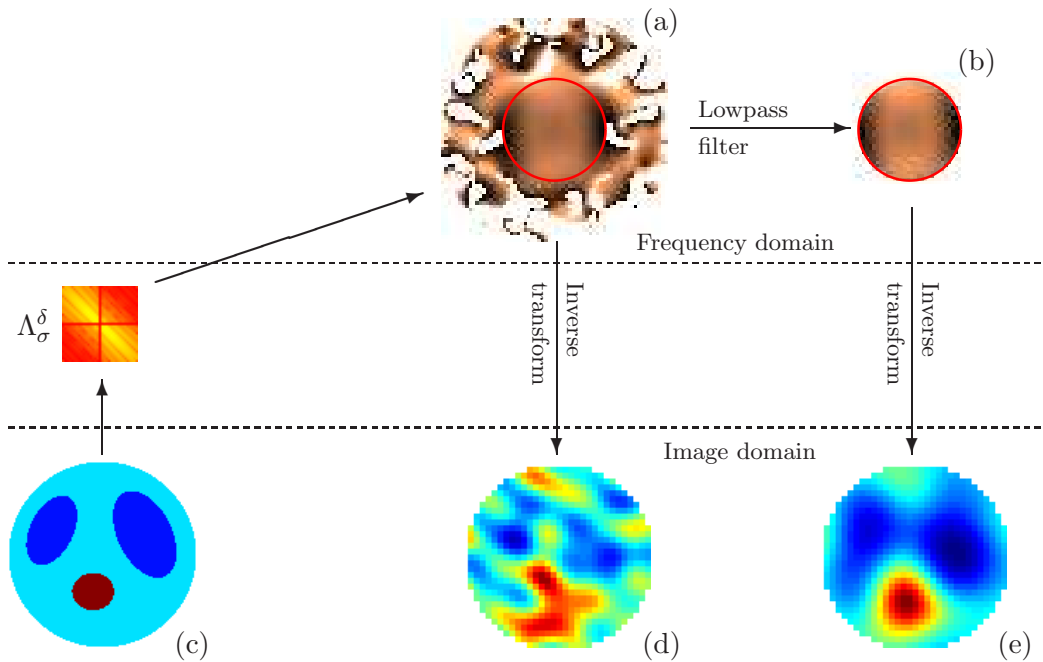
\begin{figure}[t]
\begin{picture}(320,240)
\put(-30,2){\includegraphics[height=2.5cm]{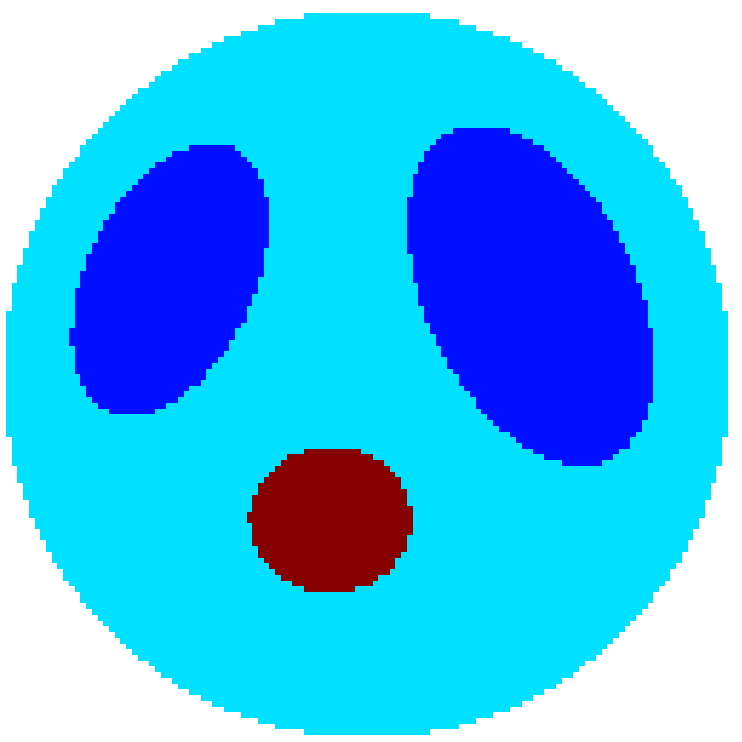}}
\put(273,2){\includegraphics[height=2.5cm]{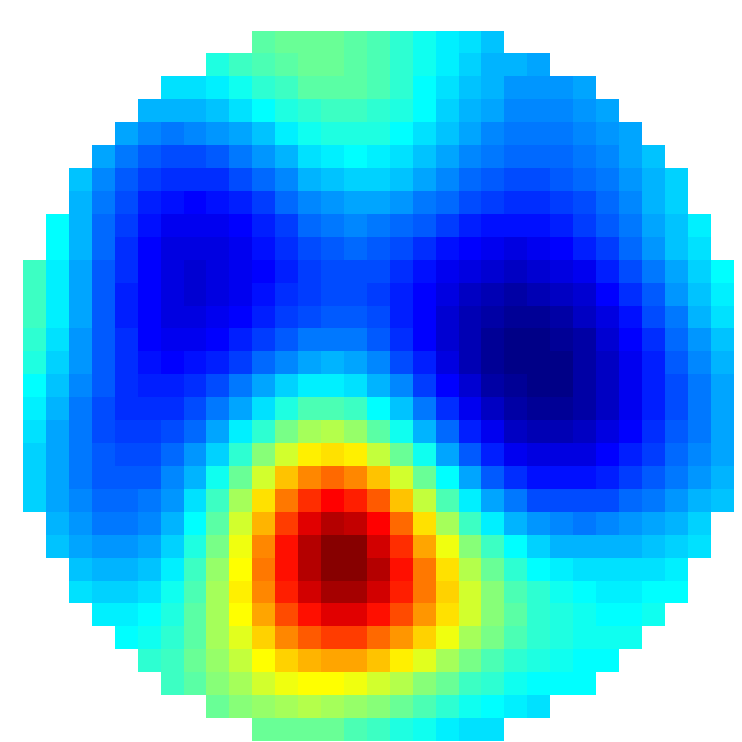}}
\put(140,2){\includegraphics[height=2.5cm]{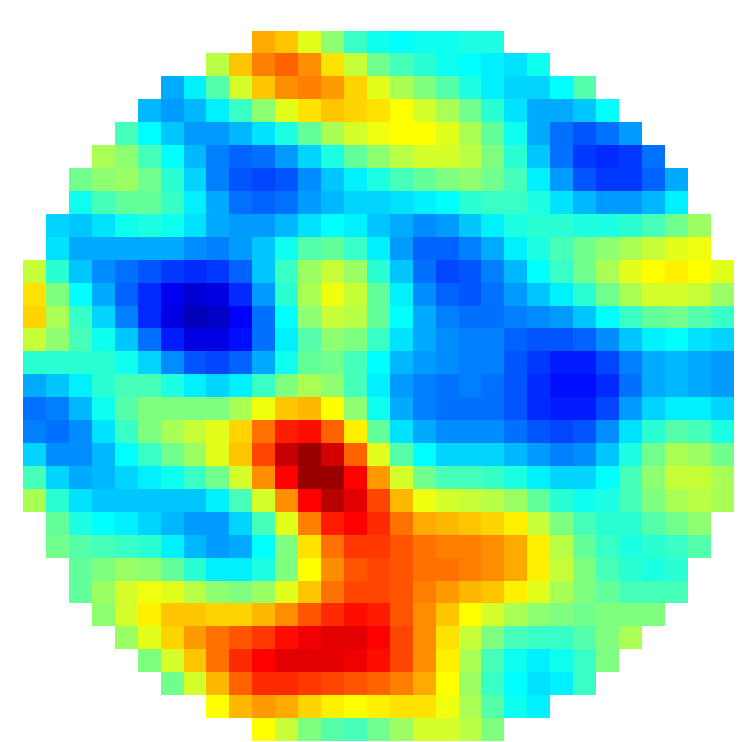}}
\put(132,155){\includegraphics[height=3cm]{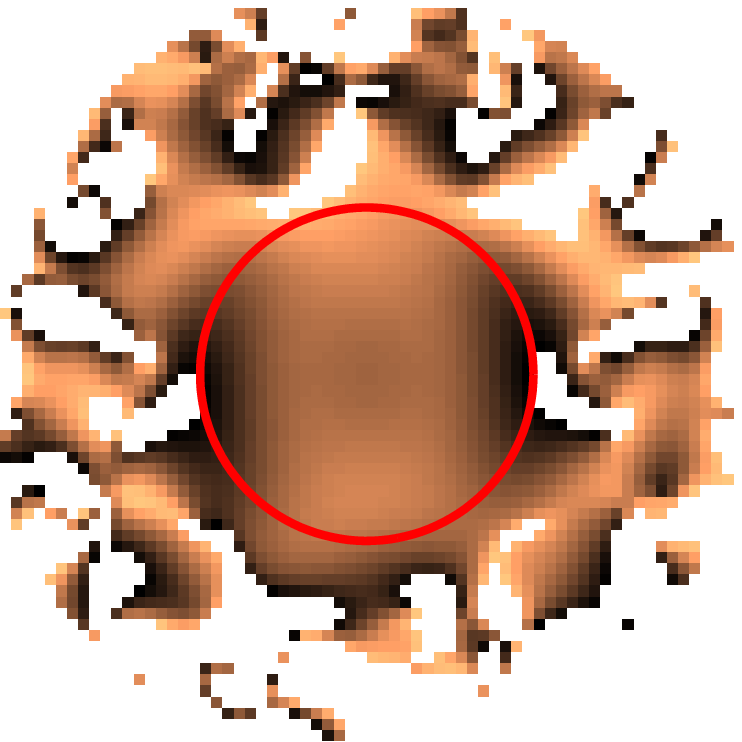}}
\put(265,155){\includegraphics[height=3cm]{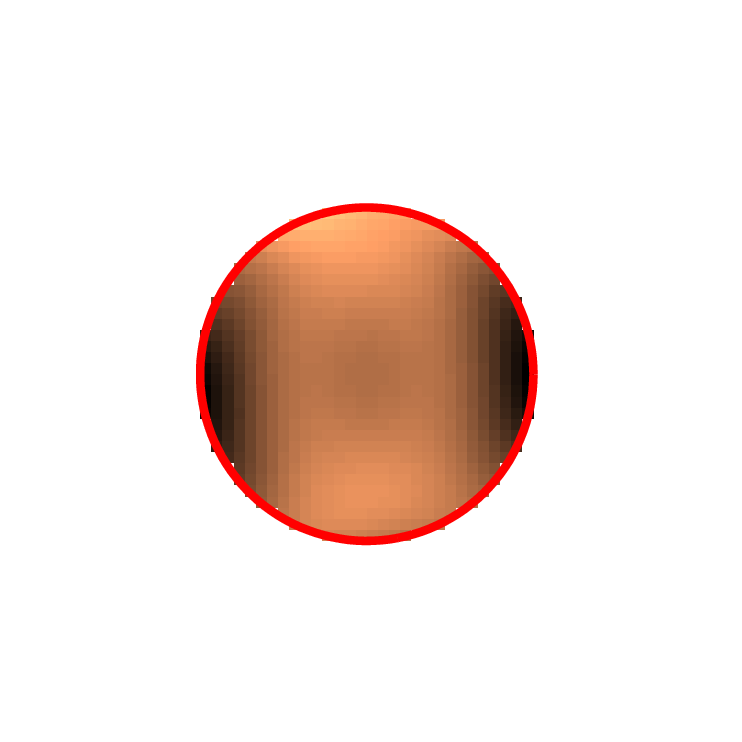}}
\put(-10,110){\includegraphics[height=1cm]{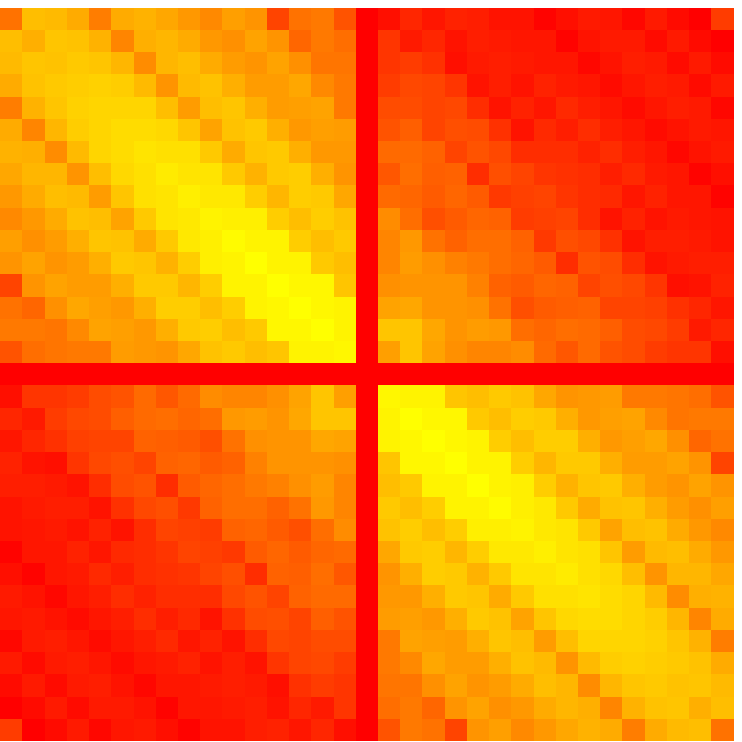}}
\thicklines
\linethickness{.1mm}
\put(308,168){\vector(0,-1){93}}
\put(312,140){\rotatebox{-90}{\footnotesize Inverse}}
\put(299,140){\rotatebox{-90}{\footnotesize transform}}
\put(175,155){\vector(0,-1){80}}
\put(179,140){\rotatebox{-90}{\footnotesize Inverse}}
\put(166,140){\rotatebox{-90}{\footnotesize transform}}
\put(5,75){\vector(0,1){28}}
\put(25,140){\vector(3,1){100}}
\put(-27,120){$\Lambda_\sigma^\delta$}
\put(225,198){\vector(1,0){60}}
\put(228,202){\footnotesize Lowpass}
\put(228,188){\footnotesize filter}
\thinlines
\multiput(-30,147)(4,0){97}{\line(1,0){2}}
\multiput(-30,85)(4,0){97}{\line(1,0){2}}
\put(215,73){\footnotesize Image domain}
\put(205,152){\footnotesize Frequency domain}
\put(207,235){(a)}
\put(325,220){(b)}
\put(32,0){(c)}
\put(207,0){(d)}
\put(340,0){(e)}
\end{picture}
\caption{\label{fig:scheme}Schematic illustration of the nonlinear low-pass filtering approach to regularized 2D EIT. The simulated heart-and-lungs phantom (c) gives rise to a finite voltage-to-current matrix $\Lambda_{\sigma}^\delta$ (orange square), which can be used to determine the nonlinear Fourier transform (a). Measurement noise causes numerical instabilities in the transform (see the irregular white patches in (a)), leading to an unstable and inaccurate reconstruction (d). However, multiplying the transform by the characteristic function of the disc $|k|<R$ yields a lowpass-filtered transform (b), which in turn gives a noise-robust approximate reconstruction (e).}
\end{figure}

There exists a nonlinear Fourier transform $\T:\C\rightarrow\C$ that is intimately connected to EIT. Namely, Nachman showed in \cite{Nachman1996} that one can use infinite-precision EIT data $\Lambda_\sigma$ to completely determine $\T$, and then apply the inverse transform, via solving a D-bar equation in the scattering variable, to recover the conductivity. However, in practice one never has such infinite-precision noise-free data; instead one works with noisy data $\Lambda_\sigma^\delta$.  The basic structure of the regularized D-bar method is shown in Figure \ref{fig:scheme}. Practical data $\Lambda_\sigma^\delta$ only allow stable computation of the nonlinear Fourier transform in a disc centered at the origin in the frequency domain. One can then use the good part of the transform in the inversion, corresponding to a nonlinear low-pass filtering. The cut-off frequency $R$ of the nonlinear low-pass filter depends logarithmically on the noise amplitude, tending to infinity asymptotically in the zero-noise limit. Analogously to linear low-pass filtering, the resulting image is smoothed and appears blurred.


\begin{figure}[t!]
\begin{picture}(360,160)
\put(0,131){\em Iteration}
\put(0,120){\em number}
\put(112,125){0}
\put(153,125){15}
\put(198,125){30}
\put(243,125){45}
\put(288,125){60}
\put(333,125){75}

\put(0,91){\em Segmentation}
\put(0,80){\em by AT flow}
\put(95,75){\includegraphics[height=1.5cm]{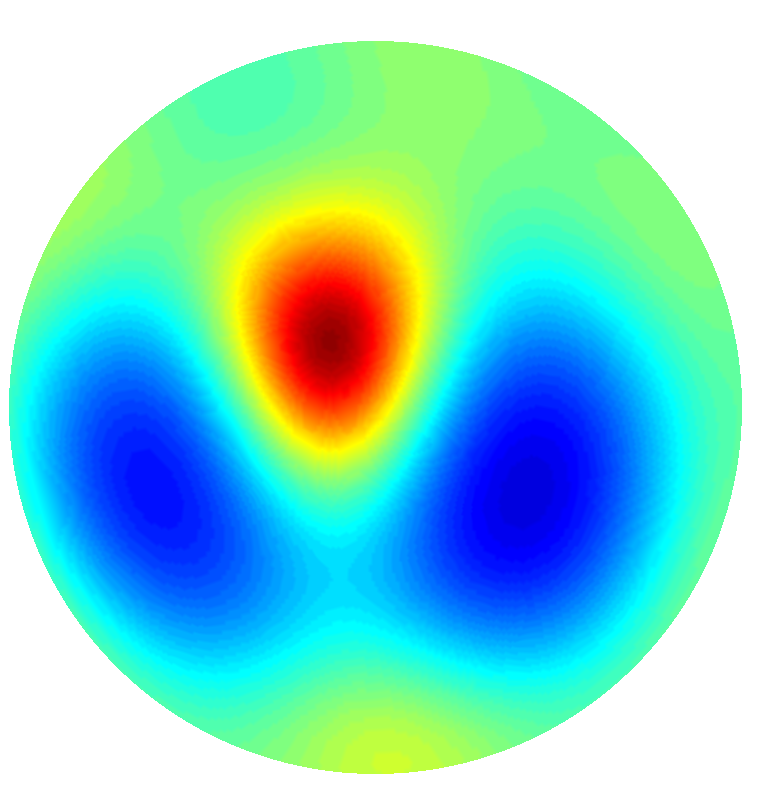}}
\put(140,75){\includegraphics[height=1.5cm]{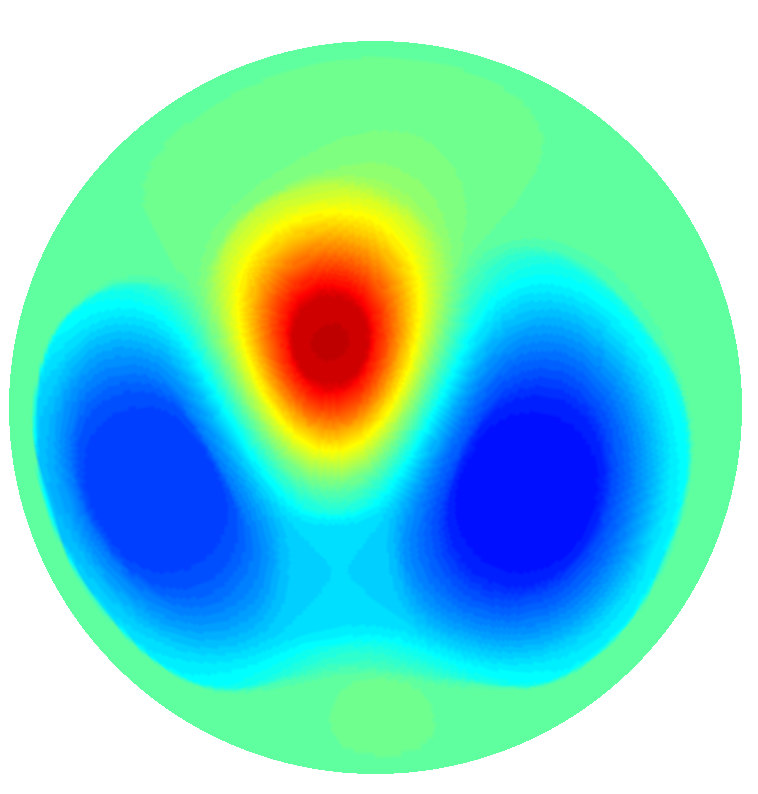}}
\put(185,75){\includegraphics[height=1.5cm]{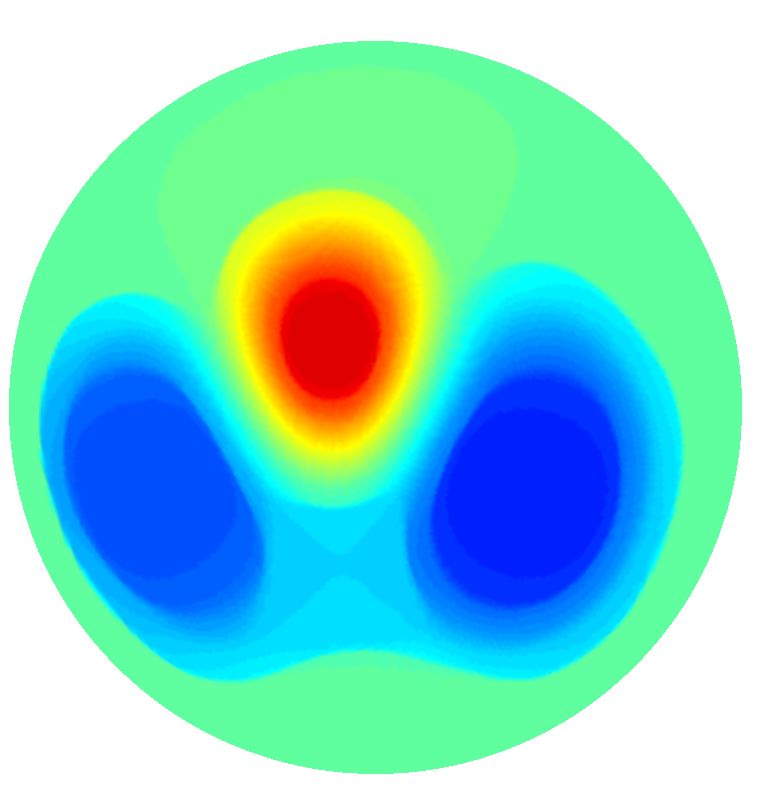}}
\put(230,75){\includegraphics[height=1.5cm]{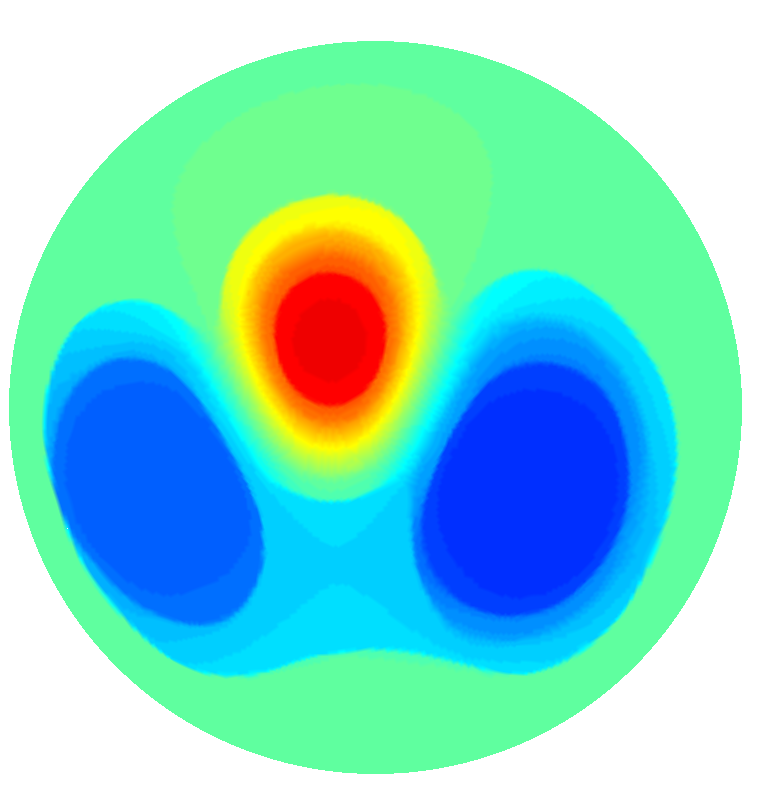}}
\put(275,75){\includegraphics[height=1.5cm]{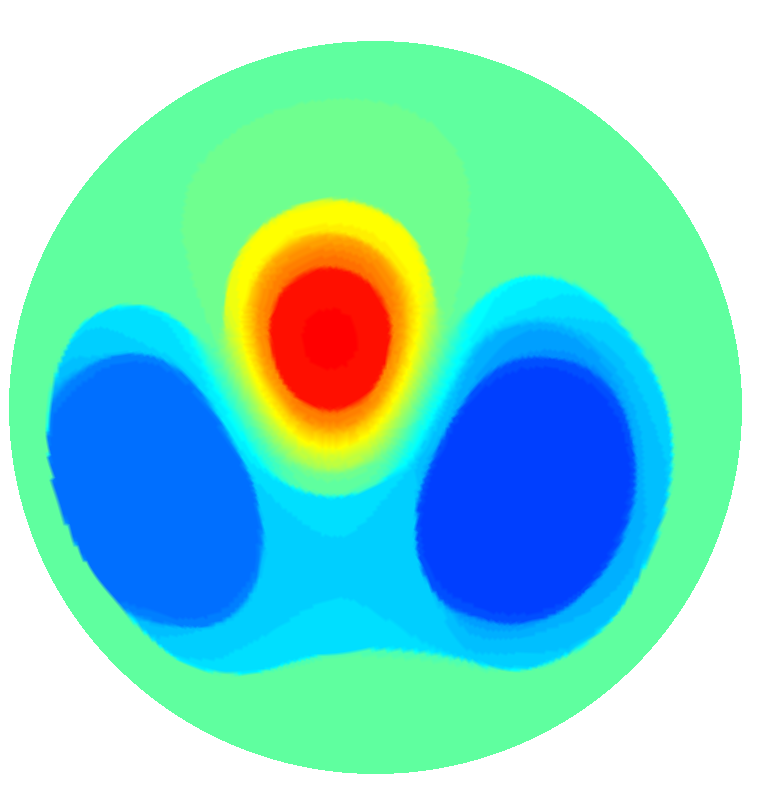}}
\put(320,75){\includegraphics[height=1.5cm]{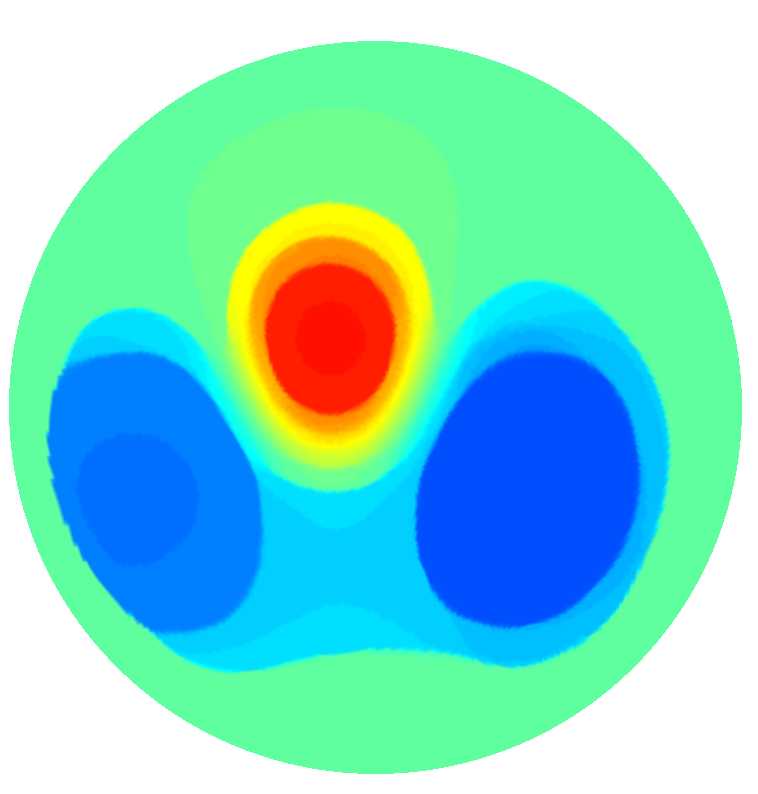}}
%
\put(0,51){\em Enhanced}
\put(0,40){\em contrast}
\put(95,30){\includegraphics[height=1.5cm]{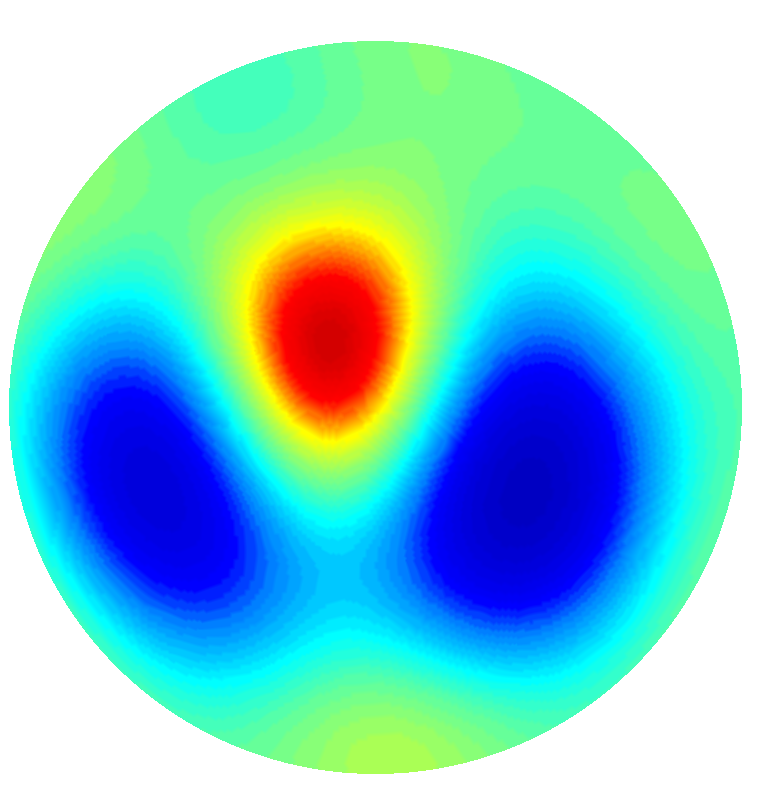}}
\put(140,30){\includegraphics[height=1.5cm]{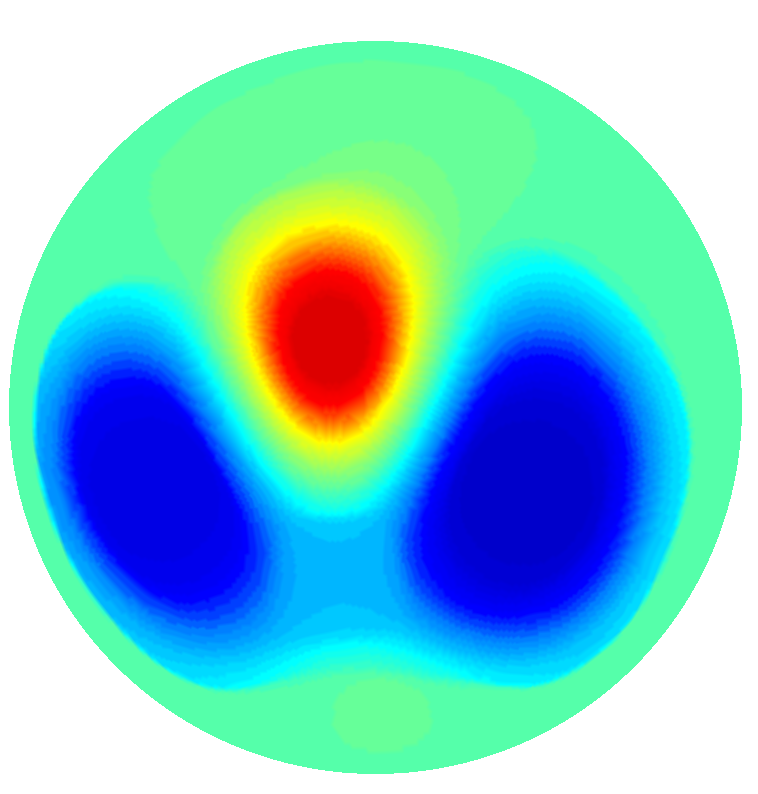}}
\put(185,30){\includegraphics[height=1.5cm]{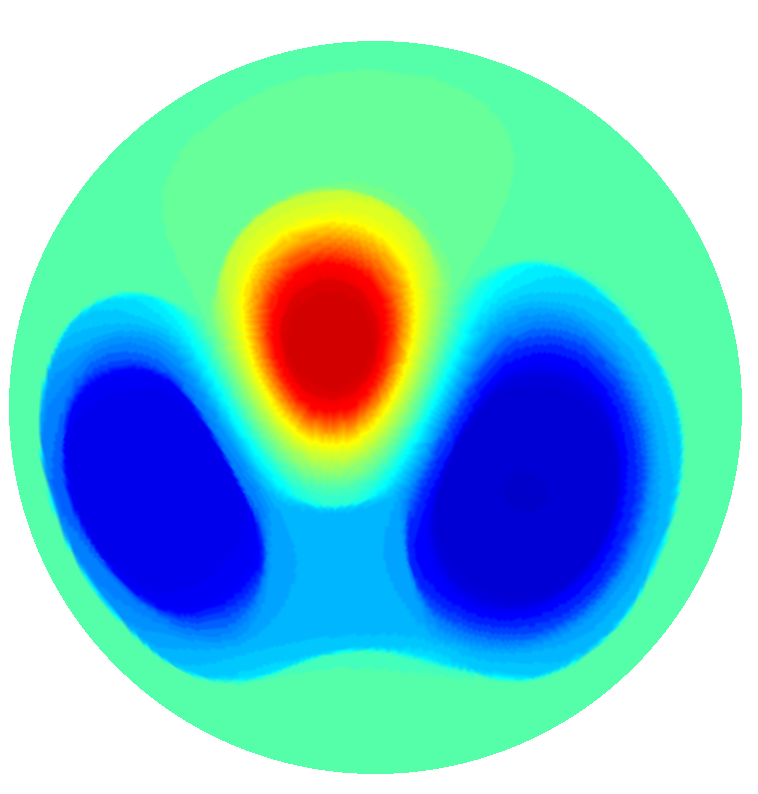}}
\put(230,30){\includegraphics[height=1.5cm]{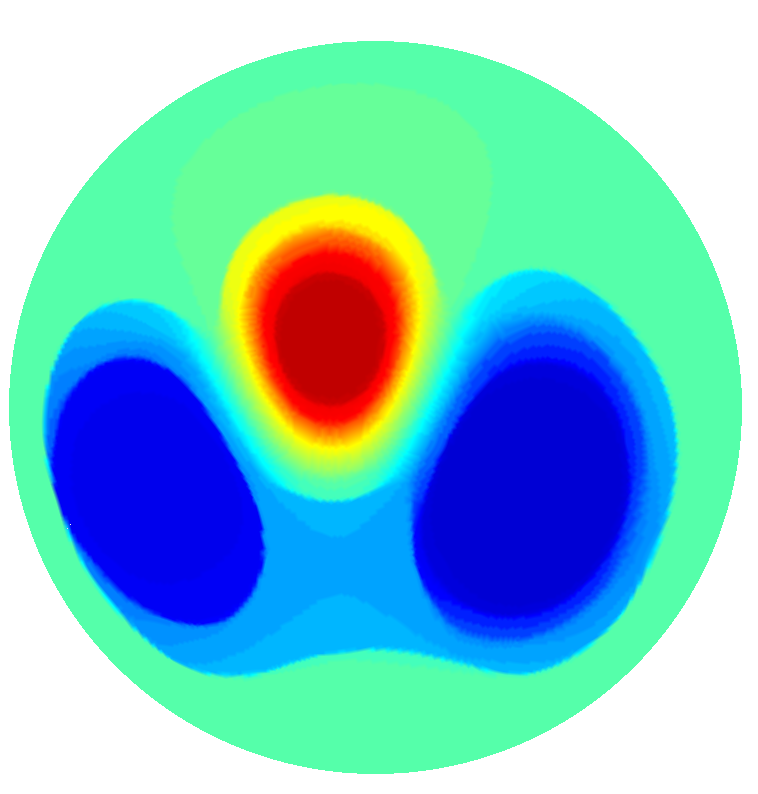}}
\put(275,30){\includegraphics[height=1.5cm]{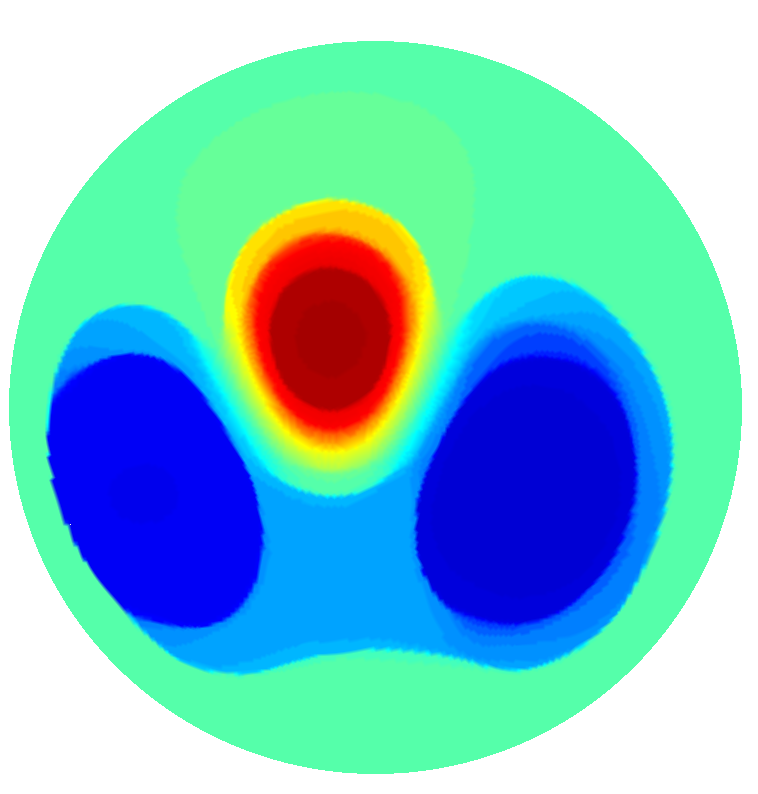}}
\put(320,30){\includegraphics[height=1.5cm]{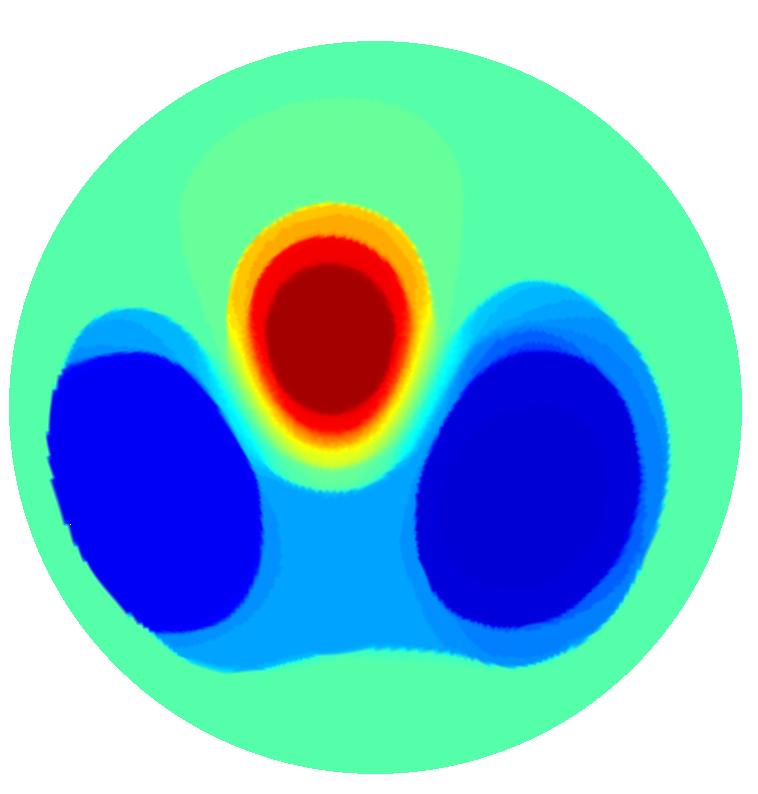}}
%
\put(0,11){\em Error in {\sc CGO}}
\put(0,0){\em sinogram (\%)}
\put(105,5){20.3}
\put(148,5){17.76}
\put(193,5){17.41}
\put(238,5){17.28}
\put(283,5){17.36}
\put(328,5){17.56}

\thinlines
\put(227,1){\line(0,1){72}}
\put(227,1){\line(1,0){46}}
\put(273,1){\line(0,1){72}}
\put(227,73){\line(1,0){46}}
\end{picture}
\caption{\label{fig:algorithmoverview}Graphical overview of the proposed EIT reconstruction algorithm. The starting point of the Ambrosio-Tortorelli segmentation flow is the outcome of the low-pass filtered D-bar method (iteration 0), here resulting from simulated EIT measurements with \newtext{0.5\% noise}. The final reconstruction is chosen to be the contrast-enhanced image having the smallest {\sc CGO} sinogram error. The measured EIT data is used in the calculation of the error.}
\end{figure}

In general, noise-robust solutions of ill-posed inverse problems are based on complementing indirect and unstable measurement data with {\em a priori} information. The transform domain of the D-bar method organizes the measurement information neatly into a stable part ($|k|<R$) and unknown/unstable part ($|k|\geq R$). Roughly speaking, information about the high frequencies of the unknown conductivity are missing from EIT data. 

What kind of {\em a priori} information do we have? In many applications of EIT it is reasonable to assume that the conductivity is piecewise constant with crisp boundaries between the homogeneous regions. Those boundaries have significant high-frequency content which is not stably represented in the data. Modern imaging science provides several options for sharpening blurred images, which may help to recover the edges in the true conductivity.  The most straightforward edge-sharpening approach, the Perona-Malik anisotropic diffusion approach \cite{Perona1990} (see \cite{Shah1996,Weickert1998} for its generalizations) proved insufficient in our work due to either too smooth edges or instability at high gradients. We therefore used the technique proposed by Ambrosio and Tortorelli \cite{Ambrosio1990} as an approximation to the classical Mumford-Shah image segmentation functional \cite{Mumford1985}.  \newtext{With this method, a}reas separate nicely and develop constant values, while the gradient can be controlled by the auxiliary variable of the model. This functional has been widely used in image processing, see for instance \cite{Chambolle1995a,Erdem2009,Jung2011}. On the other hand the approach is rather uncommon for inverse problems, there are a few examples in signal restoration \cite{Helin2011}, X-ray tomography \cite{Ramlau2007}, and particular for EIT imaging \cite{Rondi2001}. \newtext{A recent paper utilizes a Perona-Malik prior for Bayesian inversion in EIT imaging \cite{Harhanen2014}.}

It is well-known that minimizing the Ambrosio-Tortorelli (AT) functional  can be used to sharpen and segment blurry images \cite{Ambrosio1990,Shah1996}. The image flows in time according to a nonlinearly deforming diffusivity. High gradients in the initial image will develop sharp edges, while slowly varying regions will tend to constant-value areas in the final image. See Section \ref{sec:AD} for more information about AT.  Simply minimizing the AT functional for the blurry D-bar reconstruction will introduce edges, but we need to ensure that the change is for the better. In particular in medical imaging applications we need to avoid introducing artifacts. We propose controlling the AT flow via the measured EIT data. 

But how should one check the compatibility of the evolving conductivity to the measured data? One option would be to simulate the EIT measurement matrix $\Lambda_t$ at each time step for the evolving AT image and ensure that the distance to the measured matrix $\Lambda_\sigma^\delta$ decreases.   However, more robust control is provided by  a novel concept called the {\em CGO sinogram}. It consists of the traces of the complex geometric optics ({\sc CGO}) solutions (of the D-bar method) at the boundary, corresponding to low-frequency spectral parameters only. The {\em CGO sinogram} is stable to compute from EIT data \newtext{.  Furthermore, it appears to contain geometric information about the conductivity in a far more explicit form than the DN map (see Figure \ref{fig:CGOsinogramframes}), at least for some of the most traditionally used current patterns, e.g., trigonometric current patterns.}

While the AT flow will introduce edges, as desired, it also lowers the contrast in the image. To overcome this problem, we introduce a contrast-enhancement step. See Figure \ref{fig:algorithmoverview} for the result.  We remark that our method provides a novel connection between the PDE-based inverse problems community and the PDE-based image processing community having a potentially strong impact on both.

The remainder of the paper is organized as follows.  Sections~\ref{sec:Dbar}-\ref{sec:CGOsino} contain the theory behind the key pieces in the new edge preserving reconstruction method guided by the data-driven contrast enhancement of Section~\ref{sec:CE}.  In Section~\ref{sec:Dbar}, the D-bar method is reviewed.  The Ambrosio-Tortorelli functional used for recovering edges is described in Section~\ref{sec:AD}.  In Section~\ref{sec:CGOsino}, the novel {\em CGO sinogram} is introduced and the data-driven contrast enhancement is presented in Section~\ref{sec:CE}.  For the reader's convenience, Section~\ref{sec:alg} is dedicated to an explicit description of the algorithm.  In Section~\ref{sec:results}, the proposed algorithm is demonstrated on simulated noisy EIT measurements.  A discussion of the results is given in Section~\ref{sec:discus} and the take home message and conclusions of the paper are given in Section~\ref{sec:conclude}.

\section{A Brief Review of the D-bar Method}\label{sec:Dbar}

\noindent
By the {\em D-bar method} we refer to the EIT algorithm based on the theory introduced in \cite{Nachman1996}, first implemented in \cite{Siltanen2000} and equipped with an explicit regularization step in \cite{Knudsen2009}.
Alternative D-bar methods have since emerged which handle complex coefficients \cite{Francini2000,Hamilton2012,HM12_NonCirc}, less regular conductivities \cite{Brown1996,Knudsen2003,Knudsen2004a} and merely bounded $L^\infty$ conductivities \cite{Astala2006,Astala2006a,Astala2011}.  However, for the purpose of this article, we proceed with the well-established setting of \cite{Knudsen2009}, the only 2D approach with a proven regularization strategy.

The core idea in the D-bar method is to use a nonlinear Fourier transform tailor-made for EIT. To define the transform we need modified exponential functions, also called \emph{Complex Geometrical Optics} ({\sc CGO}) solutions.

Assume (for now) that $\sigma\in C^2(\overline{\Omega})$ and that $\sigma=1$ in a neighborhood of the boundary. \newtext{The constant non-unitary condition near $\bndry$ can be dealt with by scaling as discussed below.}  The conductivity equation \eqref{eq:cond_eq} can then be transformed, using the change of variables $\psi=\sqrt{\sigma} u$, to the Schr\"odinger equation 
\begin{equation*}
[-\Delta + q]\psi = 0.
\end{equation*}
Here we define the potential $q$ by extending $\sigma$ from $\Omega$ to all of $\C$ by setting $\sigma(z)\equiv 1$ for $z\in\C\setminus\Omega$ and writing
\[q(z)=\frac{\Delta\sqrt{\sigma(z)}}{\sqrt{\sigma(z)}}.\]
 Note that $q$ has compact support in $\Omega$. 

We introduce an auxiliary variable $k\in\C\setminus 0$ and look for CGO solutions  $\psi(z,k)$ satisfying
\begin{equation*}
[-\Delta + q(\cdot)]\psi(\cdot,k) = 0
\end{equation*} 
with the asymptotic condition
\[e^{-ikz}\psi(z,k)-1\in W^{1,p}(\R^2).\]
We associate $\R^2$ with $\C$ by the mapping $z=(x,y)\mapsto x+iy$, so that\linebreak $kz=(k_1+ik_2)(x+iy)$ denotes complex multiplication. For later use we introduce the related CGO solutions $ \mu(z,k)=e^{-ikz}\psi(z,k)$.

CGO solutions were originally introduced by Faddeev in \cite{Faddeev1966} and later reinvented in the context of inverse problems in \cite{Sylvester1987}. By \cite{Nachman1996} we know that the solutions $\psi$ exist and are unique for any $2<p<\infty$ for the potentials we consider here. (Other kinds of potentials may have {\em exceptional points}, or $k$ values with no unique $\psi(\,\cdot\,,k)$. See \cite{Music2013} for more details.) 

The regularized D-bar reconstruction algorithm is comprised of the following steps:
\[\Lambda^\delta_\sigma \overset{1}{\longrightarrow} \T^R(k)\overset{2}{\longrightarrow} \sigma^R(z).\]

\vspace{1em}
\begin{itemize}
\item[{\bf Step 1:}] {\bf From boundary measurements $\Lambda^\delta_\sigma$ to scattering data $\T^R$.}\\ For each fixed $k\in\C\setminus0$, solve the Fredholm integral equation of the second kind for $z\in\bndry$
\begin{equation}\label{eq:Psi_BIE}
\psi^\delta(z,k)=e^{ikz}-\int_{\bndry} G_k(z-\zeta)\left[\Lambda_\sigma^\delta - \Lambda_1\right] \psi^\delta(\zeta,k)\; dS(\zeta), 
\end{equation}
where $\Lambda_1$ is the DN map for the constant unit conductivity,  and $G_k$ is the Faddeev Green's function \cite{Faddeev1966}, with asymptotics matching $\psi$, defined by 
\[G_k(z):=e^{ikz}g_k(z),\qquad g_k(z):= \frac{1}{(2\pi)^2}\int_{\R^2} \frac{e^{iz\cdot\xi}}{|\xi|^2 +2k\left(\xi_1+i\xi_2\right)}\; d\xi.\]

\vspace{0.5em}

\noindent Evaluate the scattering transform $\T^R$ for the cut-off frequency $R>0$ using the boundary traces $\psi^\delta$ from \eqref{eq:Psi_BIE}
\begin{equation}\label{eq:tBIE}
\T^R(k):=\begin{cases}
\int_{\bndry} e^{i\bar{k}\bar{z}}\left[\Lambda_\sigma^\delta - \Lambda_1\right] \psi^\delta(z,k)\;dS(z), & |k|<R\\
0 & |k|\geq R.
\end{cases}
\end{equation}

\vspace{2em}
\item[{\bf Step 2:}] {\bf From scattering data $\T^R$ to conductivity $\sigma^R$.}\\
Fix $z\in\Omega$ and solve the D-bar equation 
\begin{equation}\label{eq:dbark_eq}
\dbar_k \mu^R(z,k) = \frac{1}{4\pi \bar{k}} \T^R(k)e(z,-k) \overline{\mu^R(z,k)}, 
\end{equation} 
where $e(z,k):=\exp\left\{i\left(kz + \bar{k}\bar{z}\right)\right\}$, saving $\mu^R(z,0)$.  The regularized conductivity $\sigma^R$ is recovered by
\begin{equation}\label{eq:muTosig}
\sigma^R(z)=\left(\mu^R(z,0)\right)^2.
\end{equation}
\end{itemize}

If $\sigma\neq 1$ near $\bndry$, but is instead a constant $\sigma_0$, the entire problem can be scaled as follows.  Let $\tilde{\sigma}=\sigma/\sigma_0$ denote a scaled conductivity which then has a value of 1 near $\bndry$. The corresponding scaled DN map is then computed by
\[\Lambda_{\tilde{\sigma}} = \sigma_0 \Lambda_\sigma.\]
After recovering $\tilde{\sigma}^R$ from Step~2 of the D-bar algorithm above, undo the scaling by multiplying by $\sigma_0$ yielding the correct $\sigma^R$.

The algorithm described above has been used successfully on both simulated and experimental data \cite{Siltanen2000,Mueller2003,Isaacson2006,Knudsen2009}.  In many applications it is well known that high frequency features such as jump discontinuities and clear edges are present, e.g. organ boundaries.  However, the low-frequency truncation needed in the regularized D-bar algorithm results in smoothed/blurred reconstructions where these high-frequency features are often absent.  This calls for post-processing of the D-bar reconstruction to reintroduce the missing features.  We propose minimizing a functional in a process known as \emph{Ambrosio-Tortorelli} image segmentation.

\section{Diffusive Image Segmentation}\label{sec:AD}
\noindent
Consider the following image processing problem. We begin with a smooth image $\widetilde{u}$, which is the result of a blurring process applied to a clean, piecewise constant image $\mathfrak{u}$. How can we recover $\mathfrak{u}$ from $\widetilde{u}$?

In 1985, Mumford and Shah introduced the following functional for the purpose of detecting boundaries in general images \cite{Mumford1985}:
\begin{equation}\label{eqn:Mumford-Shah}
E_{MS}(u,K)= \int_{\Omega\backslash K} |\nabla u|^2 dx +\beta \int_\Omega (u-\widetilde{u})^2 dx + \alpha|K|,
\end{equation}
where $K$ denotes a curve segmenting $\Omega$, $|K|$ the length of $K$, and the two  parameters $\alpha,\beta>0$ are used for weighting the terms. The idea is to find the minimum of $E_{MS}(u,K)$ over images $u$ and curves $K$; the minimizing image is then considered an edge-preserving reconstruction of $\mathfrak{u}$ from $\widetilde{u}$.

As $K$ is unknown and singular, numerically minimizing \eqref{eqn:Mumford-Shah} is a challenging task; in particular formulating a gradient-descent method with respect to $K$ is not straightforward.  Therefore, Ambrosio and Tortorelli \cite{Ambrosio1990}  proposed an elliptic approximation to \eqref{eqn:Mumford-Shah} by introducing an edge-strength function $v:\Omega\to [0,1]$ for controlling the gradient of $u$. The Ambrosio-Tortorelli (AT) functional is defined by
\begin{equation}
\label{eqn:AT_functional}
E_{AT}(u,v)=\int_\Omega  \beta(u-\widetilde{u})^2+ v^2 |\nabla u|^2 +\alpha\left( \rho |\nabla v|^2 + \frac{(1-v)^2}{4\rho} \right) dx.
\end{equation}
The additional parameter $\rho>0$ specifies, roughly speaking, the edge width of $u$. Then $E_{AT}$ $\Gamma$-converges to $E_{MS}$ as $\rho \to 0$ \cite{Ambrosio1990}, which can be understood as solutions of \eqref{eqn:AT_functional} converge to solutions of \eqref{eqn:Mumford-Shah} with the parameter $\rho\to 0$, for further explanations see \cite{Chambolle1995a}.

The advantage of \eqref{eqn:AT_functional} over \eqref{eqn:Mumford-Shah} is that the minimizer can be obtained by an artificial time evolution formulated via a coupled PDE as the gradient-descent equations with an imposed homogeneous Neumann boundary condition:
\begin{equation}
\label{eqn:AT_flow}
\left\lbrace \begin{aligned}
\partial_t u &= \diver(v^2\nabla u) - \beta(u-\widetilde{u})\hspace{0.5 cm}\text{in } \Omega\times(0,T],\\
\partial_t v &= \rho \Delta v-\frac{v|\nabla u|^2}{\alpha}+\frac{1-v}{4\rho}\hspace{0.4 cm}\text{in } \Omega\times(0,T],\\
\partial_n u&=0,\ \partial_n v=0 \hspace{2.4 cm}\text{on }  \partial\Omega\times(0,T],\\
u(\cdot,0)&=\widetilde{u},\ v(\cdot,0)=v_0 \hspace{1.9 cm}\text{in } \Omega.
\end{aligned}\right. 
\end{equation}
The equations in \eqref{eqn:AT_flow} are referred to as the {\em AT flow}.  With these gradient descent equations, implementing the numerical minimization algorithm is now a straightforward task (use finite differences or a parabolic finite element solver).  Numerous modifications of the functionals \eqref{eqn:Mumford-Shah} and \eqref{eqn:AT_functional} have been proposed, e.g. substituting the squared norm $|\nabla u|^2$ with $|\nabla u|$ in the spirit of Total Variation \cite{Alicandro1998}, or adjusting the auxiliary function $v$ as in \cite{Erdem2009}.

The effect of \eqref{eqn:AT_functional} can be visualized as follows.  First, keeping $v$ fixed, one sees that the first equation \eqref{eqn:AT_flow} minimizes the functional
\begin{equation*}
\int_\Omega \beta (u-\widetilde{u})^2 + v^2|\nabla u|^2 dx=\beta \|u-\widetilde{u}\|_{L^2(\Omega)}^2 + \|v|\nabla u|\|_{L^2(\Omega)}^2.
\end{equation*}
The representation on the right hand side is a typical form for stating regularization problems. The first fidelity term controls the distance of the obtained solution from the initial reconstruction, whereas the regularization term weights $\nabla u$ with respect to the positive (here fixed) edge-strength function $v$. The balance between the fidelity and regularization terms can be controlled by the {\em regularization parameter} $\beta>0$.

In the same manner, keeping $u$ fixed, from the second equation in \eqref{eqn:AT_flow}, we have a minimizer for a functional that can be written as
\begin{equation*}
\int_\Omega 4\rho |\nabla v|^2 + \frac{1+\frac{4\rho}{\alpha}|\nabla u|^2}{\rho}\left(\frac{1}{1+\frac{4\rho}{\alpha}|\nabla u|^2}-v\right)^2 dx.
\end{equation*}
Here we clearly see the auxiliary variable can be interpreted as a smooth approximation to the Perona-Malik filter \cite{Perona1990}, given by
\begin{equation}
\label{eqn:Perona-Malik}
g(|\nabla u|^2)=\frac{1}{1+\frac{4\rho}{\alpha}{|\nabla u|^2}}.
\end{equation}
To state a minimizing algorithm on the equations \eqref{eqn:AT_flow}, one needs to specify an initial guess of $v_0$. This is where the interpretation \eqref{eqn:Perona-Malik} comes in handy and we can set the first approximation for the auxiliary variable as $v_0=g(|\nabla \widetilde{u}|^2)$. 

Perona and Malik based their \emph{edge-aware} smoothing on the diffusion equation
\begin{equation}
\label{eqn:AD}
\partial_t u= \diver(g(|\nabla u|^2)\nabla u).
\end{equation}
Historically, the model in \eqref{eqn:AD} has a strong effect on noise removal and the smoothing of images, by keeping edges stable for a long time in the process.  However, one downside of this particular diffusion concept is that the limiting function is a constant and hence a stopping criterion during the iteration is needed \cite{Weickert1998}.  In contrast, the AT functional in \eqref{eqn:AT_functional} converges to the Mumford-Shah segmentation functional \eqref{eqn:Mumford-Shah}, for which minimizers are known to be piecewise constant with respect to the discontinuity set $K$ \cite{DeGiorgi1989}. By using the AT functional instead, the problem of defining a proper stopping criterion is shifted to choosing the correct parameters $\alpha$ and $\beta$, and making it possible to adapt the minimization problem to the reconstruction task at hand.

\newpage
\section{The CGO Sinogram}\label{sec:CGOsino}

\noindent The central idea of this study is to complement the regularized D-bar method by applying an edge-introducing image processing algorithm to a blurred EIT reconstruction. However, when manipulating the EIT image we need to ensure that we are improving the image. This is accomplished by monitoring the resulting error in a novel way (see Figure~\ref{fig:algorithmoverview}).

\begin{figure}[h]
\centering
\begin{picture}(350,320)
\put(0,25){\includegraphics[width=325 pt]{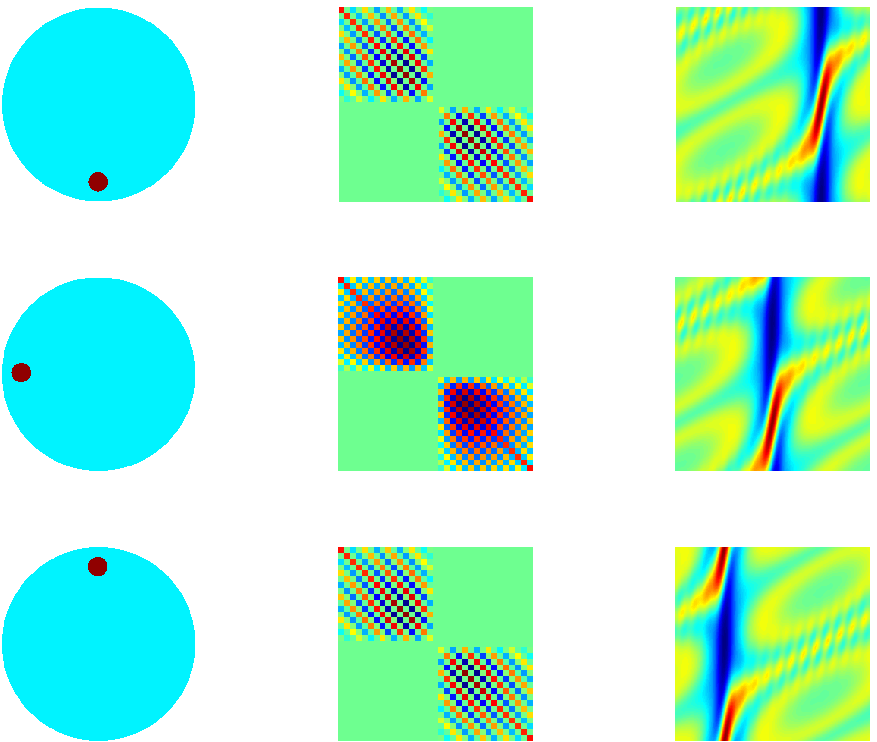}}
\put(2,305){\em Conductivity}
\put(124,305){\em DN matrix}
\put(248,305){\em CGO sinogram}
\put(300,0){\small$\mbox{arg}(z)$}
\put(328,55){\small$\mbox{arg}(k)$}
\put(26,216){\small$3\pi/2$}
\put(-6,159){\small$\pi$}
\put(28,100){\small$\pi/2$}
\put(251,226){\line(1,0){71}}
\put(251,226){\line(0,-1){3}}
\put(268.75,226){\line(0,-1){3}}
\put(286.5,226){\line(0,-1){3}}
\put(304.25,226){\line(0,-1){3}}
\put(322,226){\line(0,-1){3}}
\put(296,214){\small$3\pi/2$}
\put(251,125){\line(1,0){71}}
\put(251,125){\line(0,-1){3}}
\put(268.75,125){\line(0,-1){3}}
\put(286.5,125){\line(0,-1){3}}
\put(304.25,125){\line(0,-1){3}}
\put(322,125){\line(0,-1){3}}
\put(284,113){\small$\pi$}
\put(251,24){\line(1,0){71}}
\put(251,24){\line(0,-1){3}}
\put(268.75,24){\line(0,-1){3}}
\put(286.5,24){\line(0,-1){3}}
\put(304.25,24){\line(0,-1){3}}
\put(322,24){\line(0,-1){3}}
\put(249,12){\small$0$}
\put(262,12){\small$\pi/2$}
\put(318,12){\small$2\pi$}
\end{picture}
\caption{\label{fig:CGOsinogramframes}Left: conductivities with background one and circular inclusion of conductivity two. The polar coordinate angle of the center of the inclusion is indicated. Middle: the matrix approximation to the Dirichlet-to-Neumann map corresponding to each conductivity. More precisely, the matrix approximation to $\Lambda_\sigma^\delta-\Lambda_1$ is plotted to remove the dominating diagonal elements and bring out the differences between a homogeneous and inhomogeneous conductivity, $\sigma_1$ and $\sigma$,  more clearly. Right: CGO sinogram corresponding to each conductivity.  The color shows the values of $|\mu(z,k)-1|$ for $|z|=1$ and $|k|=2$. Note that the CGO sinograms carry clear geometric information, as opposed to the DN maps.}
\end{figure}
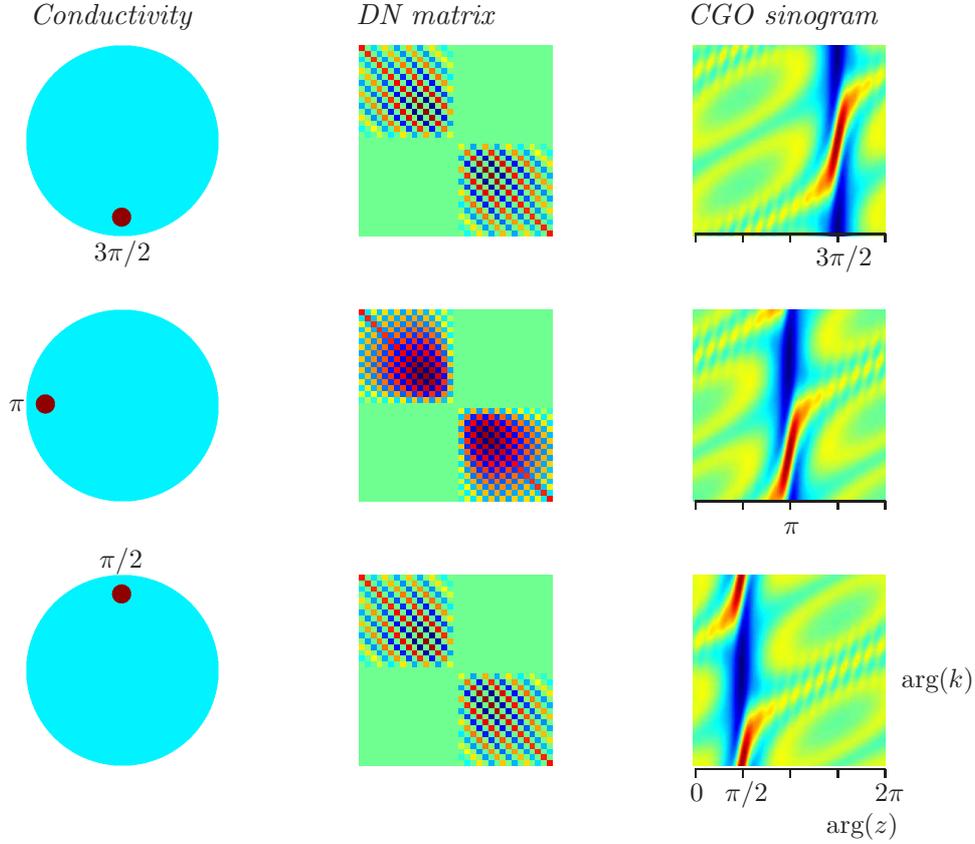

The obvious approach is to monitor the data discrepancy $\|\Lambda_\sigma^\delta-\Lambda_{\sigma^\prime}\|_Y$, where $\Lambda_\sigma^\delta$ is the measured data and $\sigma^\prime$ denotes the enhanced reconstruction, but the basic observation behind our new data fidelity term is that\newtext{, for the commonly used trigonometric current patterns,} the DN map encodes geometric information nonlinearly in a very complicated, non-intuitive, and \newtext{often} unstable way.  For example, take a look at the three conductivities shown in the left column of Figure~\ref{fig:CGOsinogramframes}. They each have one circular inclusion of conductivity two embedded into a homogeneous background of unit conductivity. The middle column shows the matrix approximations to the DN map.  Can you deduce the location of the inclusion from the DN map?  We didn't think so.

We recall from Section \ref{sec:Dbar} the related CGO solutions $\mu(z,k)=e^{-ikz}\psi(z,k)$ with asymptotic behaviour
\[\newmath{\mu(z,k)}-1\in W^{1,p}(\R^2).\]
Now take a look at the right column in Figure~\ref{fig:CGOsinogramframes}. There we show the absolute value of the difference of the {\sc CGO} solution $\mu(z,k)=\mu(e^{i\theta},2e^{i\varphi})$ and its limiting value 1, i.e. $|\mu(e^{i\theta},2e^{i\varphi})-1|$, as a function of the spectral frequency angles $\varphi$ (vertical axis) and physical angles $\theta$ (horizontal axis) where both angles range from $0$ to $2\pi$. The locations of the inclusions are immediately discernible!

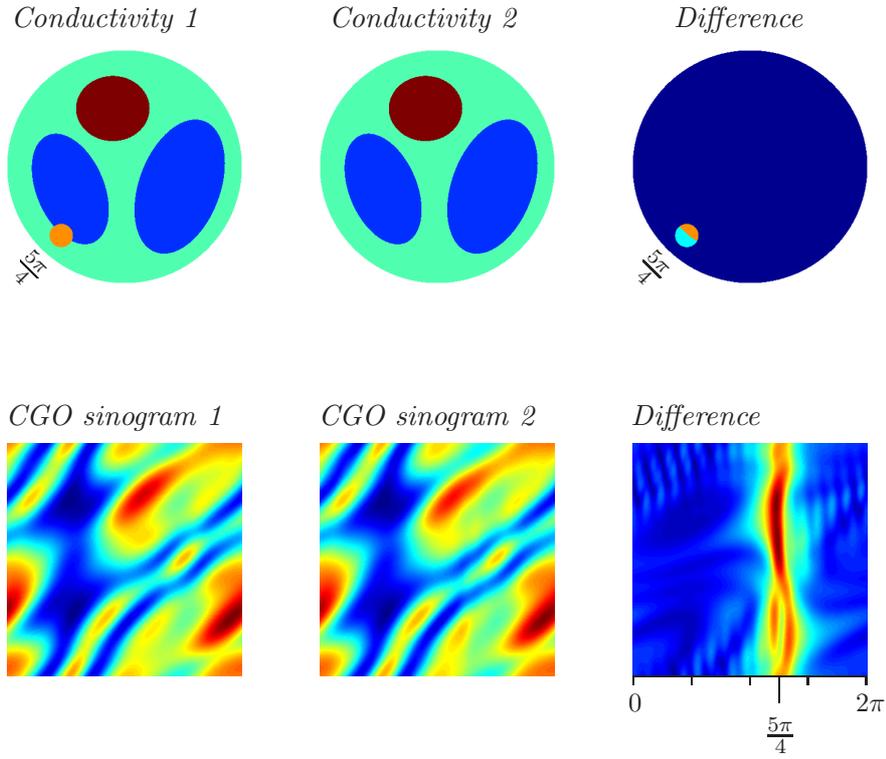
\begin{figure}[h]
\centering
\begin{picture}(350,280)
\put(0,25){\includegraphics[width=325 pt]{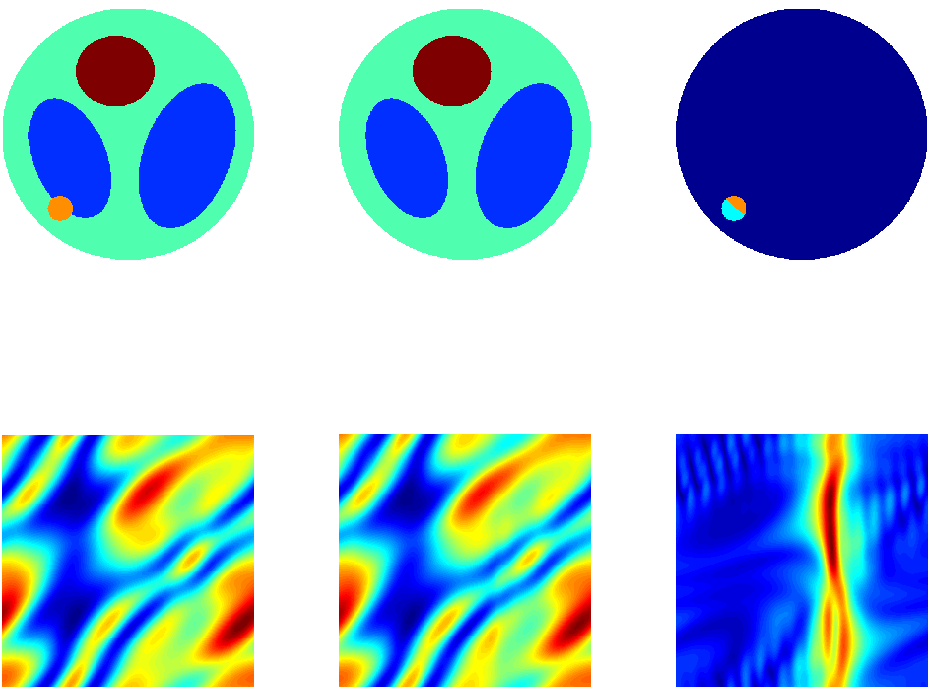}}
\put(2,270){\em Conductivity 1}
\put(121,270){\em Conductivity 2}
\put(250,270){\em Difference}
\put(233,181){\rotatebox{-45}{$\frac{5\pi}{4}$}}
\put(0,181){\rotatebox{-45}{$\frac{5\pi}{4}$}}
\put(0,120){\em CGO sinogram 1}
\put(117,120){\em CGO sinogram 2}
\put(234,120){\em Difference}
\put(234.8,25){\line(1,0){87}}
\put(256.6,25){\line(0,-1){3}}
\put(234.8,25){\line(0,-1){3}}
\put(278.4,25){\line(0,-1){3}}
\put(289.3,25){\line(0,-1){10}}
\put(300.2,25){\line(0,-1){3}}
\put(322,25){\line(0,-1){3}}
\put(233,12){\small$0$}
\put(284,0){$\frac{5\pi}{4}$}
\put(318,12){\small$2\pi$}
\end{picture}
\caption{\label{fig:CGOexample2}Top row: two heart-and-lungs conductivities (one with an extra circular inclusion), and their difference. The conductivity values are as follows. Heart: 4 (red), lungs: 1/2 (blue), background: 1 (green) and inclusion: 2 (orange). The polar coordinate angle of the center of the inclusion is indicated. Bottom row: absolute values of the corresponding CGO sinograms and their absolute difference. Note that the location of the inclusion is clearly visible in the bottom right plot.}
\end{figure}

Clearly, Figure~\ref{fig:CGOsinogramframes} is just one very simple example. See Figure~\ref{fig:CGOexample2} for a more complicated situation, where we add a circular inclusion to a heart-and-lungs phantom.  The location of the inclusion is clearly indicated in the difference of the two {\sc CGO} sinograms.

We believe that the numerical evidence presented in Figures~\ref{fig:CGOsinogramframes} and \ref{fig:CGOexample2} reflects a more general fact: calculating the {\sc CGO} sinogram 
\begin{equation}\label{def:CGOsinogram}
\begin{split}
\mathcal{S}_\sigma(\theta,\varphi,r):=&\mu(e^{i\theta},re^{i\varphi})-1 \\ 
 =&\exp(-ire^{i(\varphi+\theta)})\,\psi(e^{i\theta},re^{i\varphi})-1
\end{split}
\end{equation}
from the EIT measurements of $\sigma$ is stable (for $r>0$ inside the region of proven stability) and encodes the geometry of the conductivity in a useful and more transparent way.

\clearpage
\afterpage{\clearpage}
\section{Contrast Enhancement}\label{sec:CE}

\noindent
The Ambrosio-Tortorelli (AT) segmentation flow, discussed in Section \ref{sec:AD} above, can transform a blurry image into a sharper image. However, the AT flow comes with a reduction in contrast, which in turn is a key benefit of EIT imaging.  To overcome this obstacle, we propose using a \emph{data-driven} contrast enhancement technique kept in check by the {\sc CGO} sinogram.  

We search for a contrast enhanced (corrected) $\sigmaCE$, such that the values $\sigma(z)$ are stretched or damped according to the resulting error in CGO sinogram. We utilize a two parameter model and consider values greater and lower 1 independently, where 1 is the value near the boundary. For that let  
\[f(z):=\sigmaAT(z) -1,\]
denote the difference between the conductivity after AT flow and the constant 1.  Recall the  {\em a priori} constants $c$ and $C$ known to bound the conductivity $0 < c\leq \sigma(z)<C$.  In practice, such ballpark bounds are readily available.  For example, in chest imaging, using an applied current with frequency 100 kHz, internal conductivities range from around $0.02$ Siemens/meter (e.g. fat) to $0.71$ Siemens/meter (e.g. heart tissue) \cite{INCR}.

Denote by $m$ and $M$ the following minimum and maximum values
\[m:=\underset{z\in\Omega}{\min} f(z),\qquad  M:=\underset{z\in\Omega}{\max} f(z),\]
assuming these values are nonzero.
Define the scaling function $f_{s,t}$ for the scaling parameter $(s,t)\in [0,1]^2$ by
\begin{equation*}
f_{s,t}(z)=\begin{cases}
s\frac{f(z)}{m}(c-1) & \mbox{for $z$ satisfying $f(z)<0$}\\
t\frac{f(z)}{M}(C-1) & \mbox{for $z$ satisfying $f(z)\geq 0$}
\end{cases}
\end{equation*}
and set
\begin{equation}\label{eq:sig_st}
\sigma_{s,t}(z):=1+f_{s,t}(z).
\end{equation}

The optimal contrast enhanced conductivity $\sigmaCE$ within the bounds $c$ and $C$ is then determined by minimizing the {\sc CGO} sinogram error over the scaling parameter $(s,t)\in [0,1]^2$, i.e. 
\begin{equation}\label{eq:st_OPT}
(s_0,t_0):= \left. \argmin_{(s,t)\in[0,1]^2} \middle\{ \|\mathcal{S}^\delta_\sigma(\theta,\varphi,r)-\mathcal{S}_{\sigma_{s,t}}(\theta,\varphi,r)\|_{L^2(\mathbb{T}^2)}^2\right\},
\end{equation}
where $\mathcal{S}^\delta_\sigma(\theta,\varphi,r)$ is the CGO sinogram corresponding to the noisy measurement of $\sigma$, and $\mathcal{S}_{\sigma_{s,t}}(\theta,\varphi,r)$ is computed from $\sigma_{s,t}$. Then \eqref{eq:sig_st} is used to define
\begin{equation*}
\sigmaCE(z):=\sigma_{s_0,t_0}(z).
\end{equation*}
Various algorithms are available to find an approximation to the optimal $s_0$ and $t_0$ in (\ref{eq:st_OPT}).  In this \newtext{introductory} work we apply a derivative free pattern search algorithm which is known to be stable in the presence of noise: the {\sc DIRECT} (DIviding RECTangles) pattern search for global optimization \cite{Jones1993,Finkel2003}.

\begin{figure}[h]
\centering
\begin{picture}(300,125)
\put(0,10){\includegraphics[width=100 pt]{images/HnL_noisy_AT_it0.eps}}
\put(197,10){\includegraphics[width=108 pt]{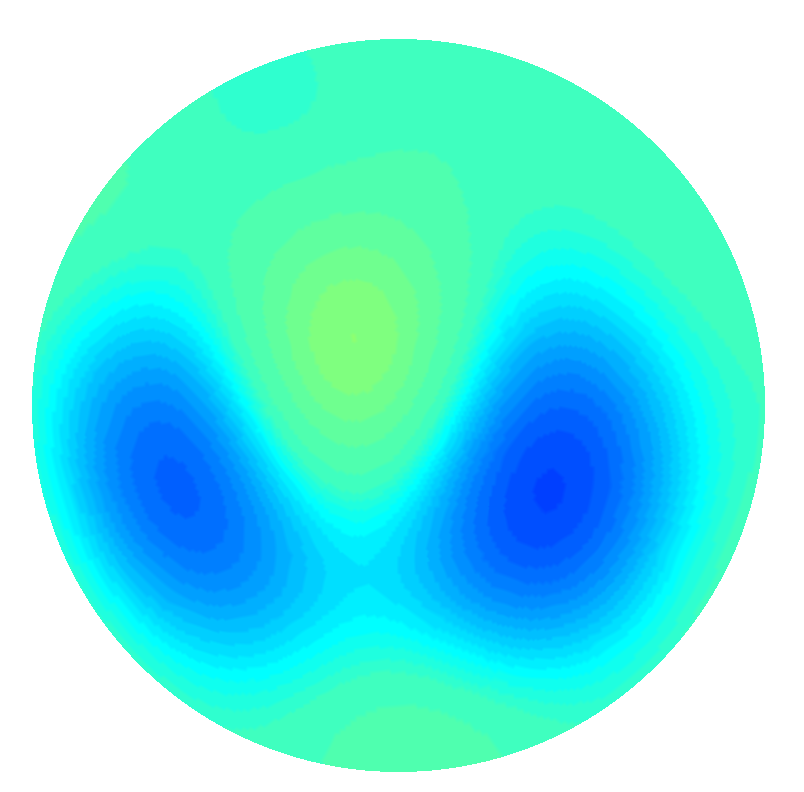}}
\put(100,10){\includegraphics[width=100 pt]{images/HnL_noisy_AT_CE_it0.eps}}
\put(32,115){\textbf{D-bar}}
\put(130,115){\textbf{CGO}}
\put(240,115){\textbf{DN}}
\put(20,0){Error 15.11\%}
\put(120,0){Error 14.72\%}
\put(220,0){Error 15.05\%}
\end{picture}
\caption{\label{fig:DN_vs_CGO}Comparison of contrast enhanced solutions \newtext{on noisy data}. Left: initial D-bar reconstruction with relative $l^1$-error 15.11\% to the original phantom. Middle: solution chosen by minimal {\sc CGO} sinogram error, relative $l^1$-error 14.72\% to the original phantom. Right: solution with minimal error in DN maps, relative $l^1$-error 15.05\% to the original phantom.}
\end{figure}

In Figure~\ref{fig:DN_vs_CGO} we see clear evidence that minimizing with respect to the {\sc CGO} sinogram is advantageous to the DN map for preserving important image features. The figure shows the heart and lungs test problem with conductivity reconstructed by the D-bar algorithm \newtext{applied to noisy EIT measurements} (left) and the optimal contrast-enhanced images chosen by minimizing the error of {\sc CGO} sinogram (middle) and DN map (right) to the measured data. Clearly, the image chosen by the {\sc CGO} sinogram minimization more accurately portrays the original than the one chosen by the DN map. Even though both solutions were able to minimize the $l^1$-error to the true phantom, from the introductory example in Figure~\ref{fig:hammer_to_the_head}, the features preserved differ immensely (i.e. the heart is barely visible in the DN guided minimization).

\section{A Data-Driven Edge-Preserving D-bar Algorithm}\label{sec:alg}
\noindent
The aim of this paper is to combine the strengths of each the methods described above in Sections~\ref{sec:Dbar}-\ref{sec:CGOsino}, with the data-driven contrast enhancement of Section~\ref{sec:CE}. We compute the reconstruction of the conductivity with the regularized D-bar method and reintroduce edges by a data-driven post-processing of the image. The post-processing is monitored by the {\sc CGO} sinogram error which incorporates geometric information about the reconstruction. We thus propose the following \textit{Data-Driven Edge-Preserving D-bar Algorithm}:\\

\begin{itemize}
\item[{\bf Step 1:}] Fix $R>0$ and compute the regularized D-bar reconstruction $\sigma^R$ using the D-bar algorithm described in Section~\ref{sec:Dbar}.

\vspace{1em}
\item[{\bf Step 2:}] Fix a radius $r$ located in the stable disc: $R>r>0$. 
\begin{enumerate}[(i)]
\item Compute the {\sc CGO} sinogram $\mathcal{S}_\sigma^\delta(\theta,\varphi,r)$ from the original noisy data $\Lambda_\sigma^\delta$ by solving the boundary integral equation (\ref{eq:Psi_BIE}) for $\psi^\delta(z,k)=e^{ikz}\mu^\delta(z,k)$ and setting
$$
\mathcal{S}_\sigma^\delta(\theta,\varphi,r)=\mu^{\delta}(e^{i\theta},re^{i\phi})-1.
$$
\item Calculate the {\sc CGO} sinogram  $\mathcal{S}_{\sigma^R}(\theta,\varphi,r)$ for the D-bar image.
\item Record the relative {\sc CGO} sinogram error
$$
E_0:=\frac{\|\mathcal{S}^\delta_\sigma(\theta,\varphi,R)-\mathcal{S}_{\sigma^R}(\theta,\varphi,R)\|_{L^2(\mathbb{T}^2)}}
{\|\mathcal{S}^\delta_\sigma(\theta,\varphi,R)\|_{L^2(\mathbb{T}^2)} }.
$$
\end{enumerate}

\vspace{1em}
\item[{\bf Step 3:}]
Reintroduce edges to $\sigma^R$ via minimization of the Ambrosio-Tortorelli functional by solving the gradient descent equations \eqref{eqn:AT_flow}. 
\begin{enumerate}[(i)]
\item Initialize the constants $\beta,\alpha$ and $\rho$. 
\item Calculate the initial approximation for the auxiliary function as \linebreak $v_0=g(|\nabla \sigma^R|^2)$ defined by \eqref{eqn:Perona-Malik}.
\item Choose a time step $t$ and begin solving \eqref{eqn:AT_flow} iteratively.
\end{enumerate}

\vspace{1em}
\item[{\bf Step 4:}] Check the flow every $J$, e.g. $J=5$, time steps as follows. For the $j$-th check:
\begin{enumerate}[(i)]
\item Denote the image to be checked $\sigmaAT_j$.
\item Determine $\sigmaCE_j$, the optimal contrast enhanced version of $\sigmaAT_j$ using the {\sc CGO} sinogram optimization method described in Section~\ref{sec:CE}.
\item Record the {\sc CGO} sinogram error $E_j$ for $\sigmaCE_j$.
\item If $E_j<E_{j-1}$, return to the AT flow with $\sigmaAT_j$ (the non-contrast enhanced version) and repeat steps (i)-(iii).  If not, or if a maximum number of iterations $J_{max}$ is reached, set 
\begin{equation*}
\sigma_{NEW}(z):=\sigmaCE_j(z)
\end{equation*}
and the algorithm is complete.
\end{enumerate}
\end{itemize}

\section{Computational results}\label{sec:results}
\noindent
We tested the algorithm on simulated noisy EIT measurement data for test cases of potential interest for the medical and industrial communities.

To simulate the EIT measurements,  we solved the Neumann problem corresponding to the conductivity equation
\begin{equation*}
\begin{array}{rcl}
\nabla \cdot\sigma\nabla u&=& 0,\quad z\in \Omega\subset\R^2\\
\sigma\frac{\partial u}{\partial\nu} &=& \phi_j, \quad z\in\bndry,
\end{array}
\end{equation*}
for $j=-16,\ldots,-1,1,\ldots,16$, representing 32 linearly independent current patterns.  In this paper $\Omega$ is the unit disc, and we used the Fourier basis functions
\[ \phi_j\left(e^{i\theta}\right)=\frac{1}{\sqrt{2}} e^{ij\theta}, \quad z=e^{i\theta}\in\bndry.\]

The matrix approximation to the DN maps $\Lambda_\sigma$ and $\Lambda_1$ (the DN map corresponding to a uniform unit conductivity) were computed for each test problem using the standard methods described in \cite{Mueller2012}. The numerical implementation of the regularized D-bar method was first fully described in \cite{Knudsen2009} and later explained in more detail in \cite[Section 13.2]{Mueller2012}, including freely available Matlab code.  

A discussion of a numerical implementation with finite differences of the AT flow can be found in \cite{Erdem2009}.  In this study, we used the Matlab integrated PDE solver for the equations \eqref{eqn:AT_flow} and for each test problem below the choice of parameters is given.  

Additional zero-mean random Gaussian noise was added to the DN matrix $\Lambda_\sigma$ so that
\begin{equation}
\label{eqn:DN_error}
\|\Lambda_\sigma-\Lambda_\sigma^\delta\|_{H^{1/2}(\bndry)\to H^{-1/2}(\bndry)}\leq \delta
\end{equation}
using the methods described in \cite{Knudsen2009}.  While types of noise and their magnitudes vary among EIT devices, a benchmark of around 0.01\% has been obtained \cite{Cheney1999}.

The resulting {\sc CGO} sinogram and conductivity errors stated below are given as relative errors. The relative error of  CGO sinograms corresponding to the original conductivity $\sigma$ and to some reconstruction $\sigma_d$ was measured by
\begin{equation}\label{sinogramerror}
\frac{\|\mathcal{S}^\delta_\sigma(\theta,\varphi,R)-\mathcal{S}_{\sigma_{d}}(\theta,\varphi,R)\|_{L^2(\mathbb{T}^2)}}
{\|\mathcal{S}^\delta_\sigma(\theta,\varphi,R)\|_{L^2(\mathbb{T}^2)} }.
\end{equation}
Similarly we measured the relative error of the reconstructed conductivities to the known original via
\begin{equation}\label{groundtrutherror}
\frac{\|\sigma-\sigma_{d}\|_{L^p(\Omega)}}
{\|\sigma\|_{L^p(\Omega)} }
\end{equation}
with $p=1$ or $p=2$. We note that \eqref{sinogramerror} can be measured for real data cases, contrary to \eqref{groundtrutherror} for which the knowledge of the correct conductivity $\sigma$ is needed.

\subsection{A Heart-and-Lungs Phantom}\label{sec:heartlungs}
\begin{figure}[h!]
\centering
\begin{picture}(300,120)
\put(0,0){\includegraphics[width=300 pt]{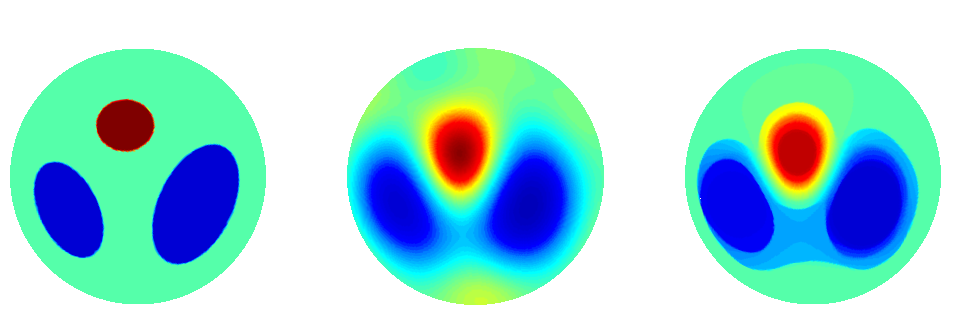}}
\put(-5,92){\textbf{True conductivity}}
\put(133,92){\textbf{D-bar}}
\put(230,92){\textbf{Improved}}
\end{picture}
\caption{\label{fig:HnL_noisy_sum}Illustration of original Heart-and-Lungs phantom with the initial D-bar reconstruction (truncation radius $R=4$) and the contrast-enhanced diffused solution of the Ambrosio-Tortorelli functional at iteration 45.}

\end{figure}

The leftmost image in Figure~\ref{fig:HnL_noisy_sum} shows our piecewise constant phantom. The lungs have conductivity 0.5, the heart has conductivity 2, and the background has unit conductivity.

The D-bar reconstruction was obtained using a truncation radius of $R=4$  in the scattering data. The D-bar method allows computing the reconstruction at arbitrary points in the $z$-plane, and we chose to reconstruct directly at the points of the FEM mesh (of 33025 elements) to be used in the AT flow.  
We have added additional noise of amplitude $\delta=0.005$ to the DN map satisfying \eqref{eqn:DN_error}, and corresponding to $0.5\%$ noise far surpassing the $0.01\%$ benchmark for measurement noise.  The D-bar reconstruction is shown in the middle image in Figure~\ref{fig:HnL_noisy_sum}; it has a relative CGO sinogram error of 20.3\%.

The AT flow (Step~3) was computed with parameters 
\[\alpha=200, \hspace{0.25cm} \beta=0.1, \hspace{0.25cm} \rho=0.1,\]
and the contrast enhanced solution $\sigmaCE$ was computed every fifth iteration using the {\sc CGO} sinogram with $r=|k|=2$, well within the \newtext{observed} reliable region for the regularized D-bar method. The boundary constants were roughly chosen as $c=0.1$ and $C=4$. A summary of the obtained solutions is illustrated in Figure \ref{fig:HnL_AT_illu}.  

The optimal solution was obtained at iteration 45 with a relative error in CGO sinogram of 17.28\%.   Figure~\ref{fig:HnL_convergence_plots} displays the error in {\sc CGO} sinogram as well as the error of reconstruction to the true conductivity throughout the evolution of the AT flow. 

\begin{figure}[h!]
\centering
\begin{picture}(310,220)
\put(10,0){\includegraphics[width=100 pt]{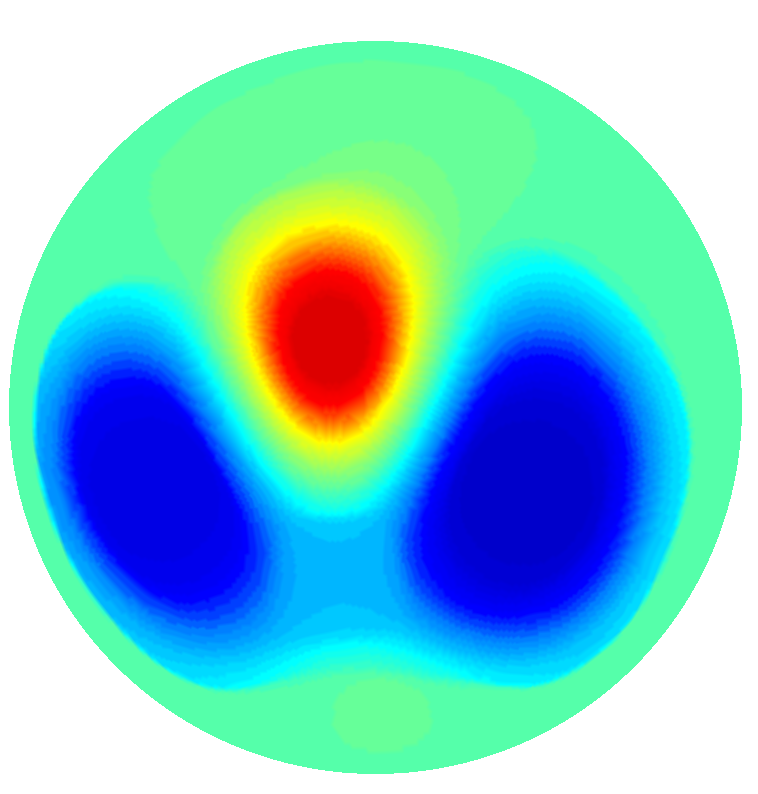}}
\put(110,0){\includegraphics[width=100 pt]{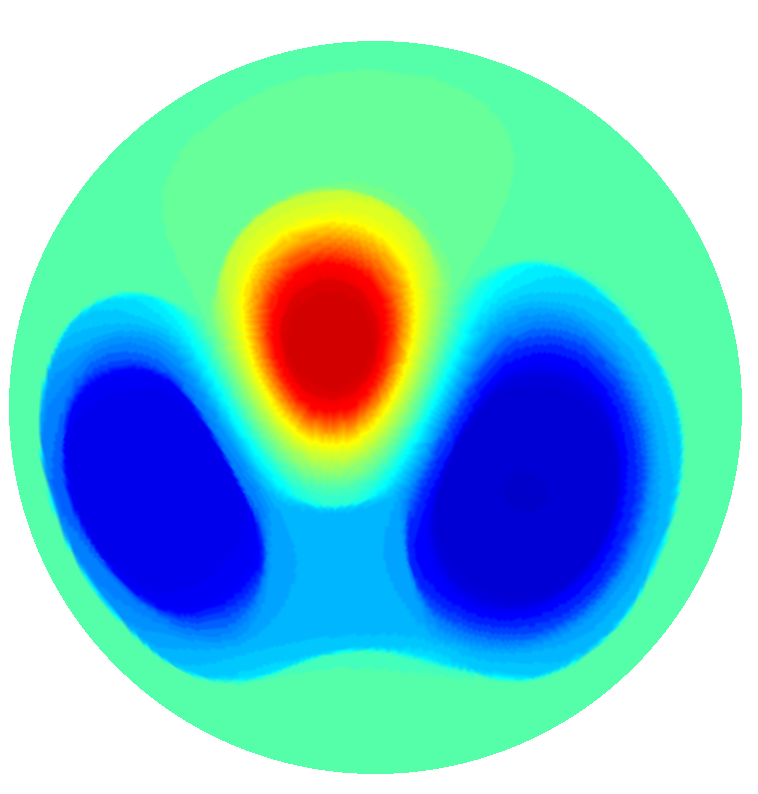}}
\put(210,0){\includegraphics[width=100 pt]{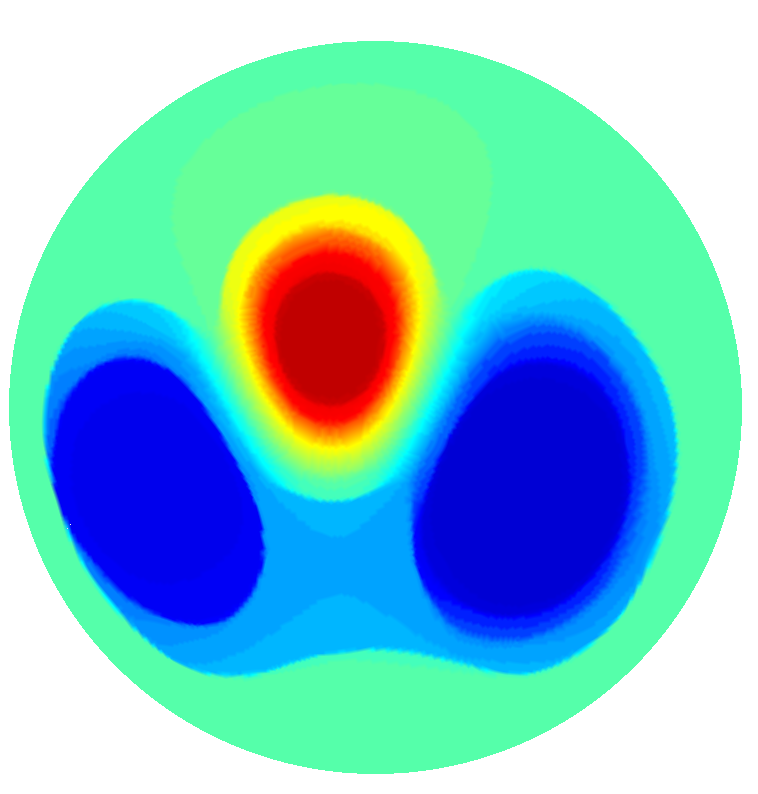}}
\put(10,100){\includegraphics[width=100 pt]{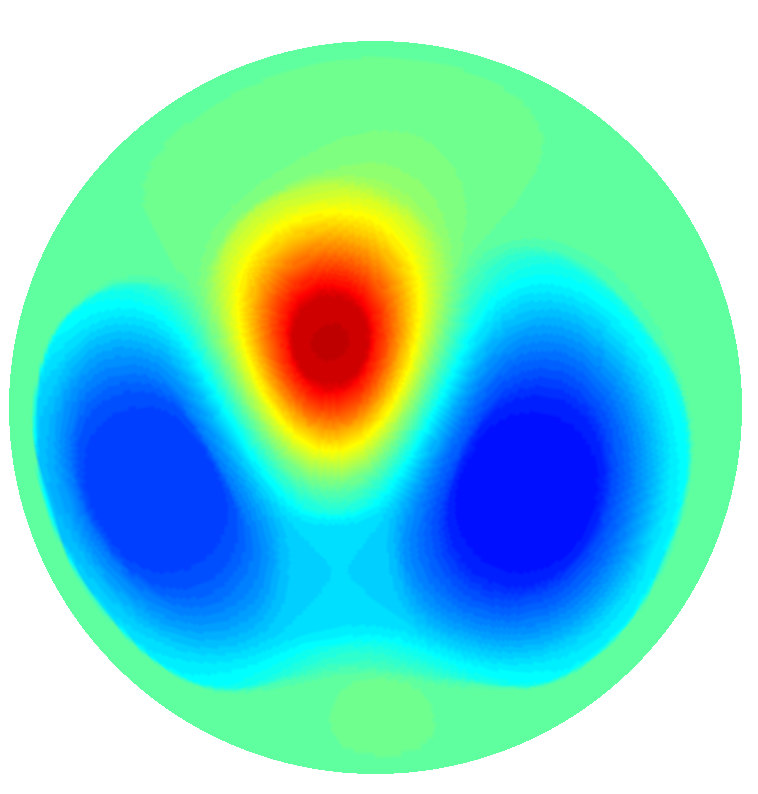}}
\put(110,100){\includegraphics[width=100 pt]{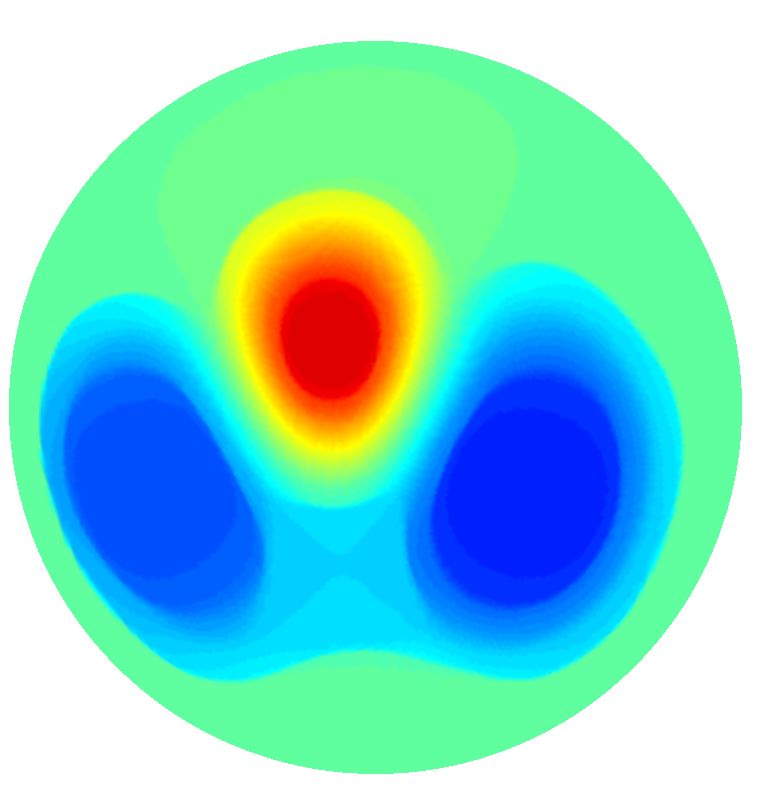}}
\put(210,100){\includegraphics[width=100 pt]{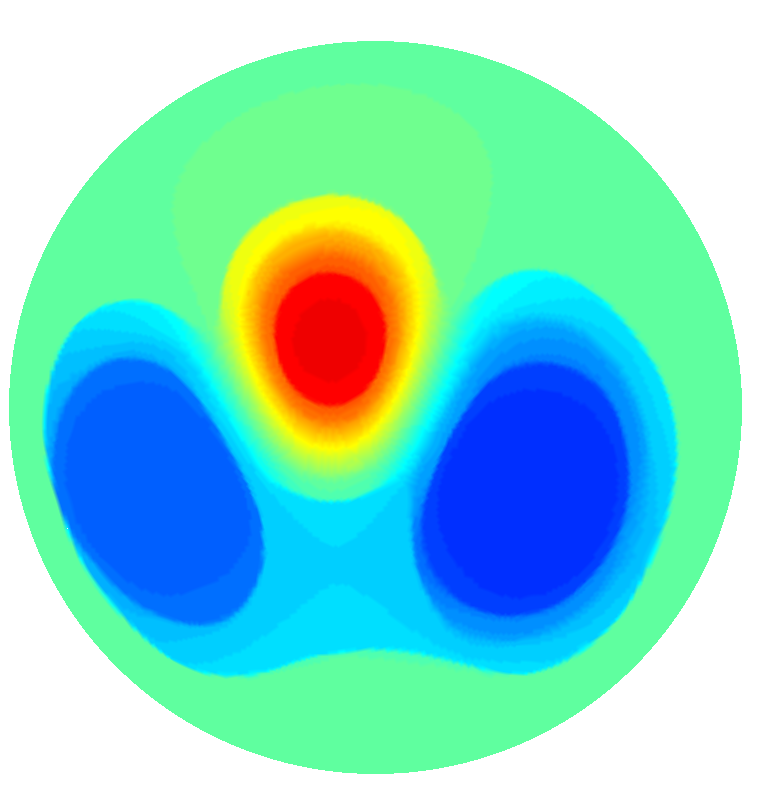}}
\put(35,205){{\scriptsize \textbf{Iteration 15}}}
\put(135,205){{\scriptsize \textbf{Iteration 30}}}
\put(235,205){{\scriptsize \textbf{Iteration 45}}}
\put(0,135){\rotatebox{90}{\footnotesize \textbf{AT flow}}}
\put(0,20){\rotatebox{90}{\footnotesize \textbf{CE solutions}}}
\end{picture}
\caption{\label{fig:HnL_AT_illu}Illustration of the minimization process of the Ambrosio-Tortorelli functional for three different stages, including the contrast enhanced solutions with smallest error in CGO sinogram.}

\end{figure} 
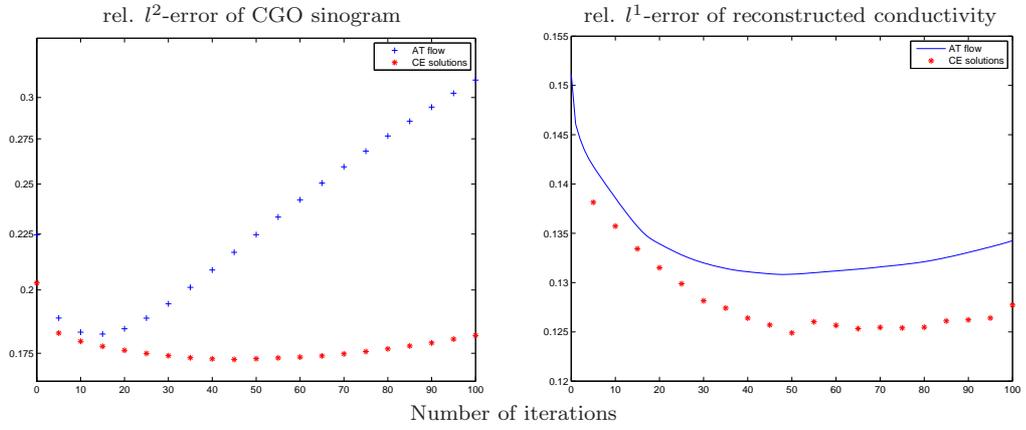
\begin{figure}[h!]
\centering
\begin{picture}(400,180)
\put(0,10){\includegraphics[width=179 pt]{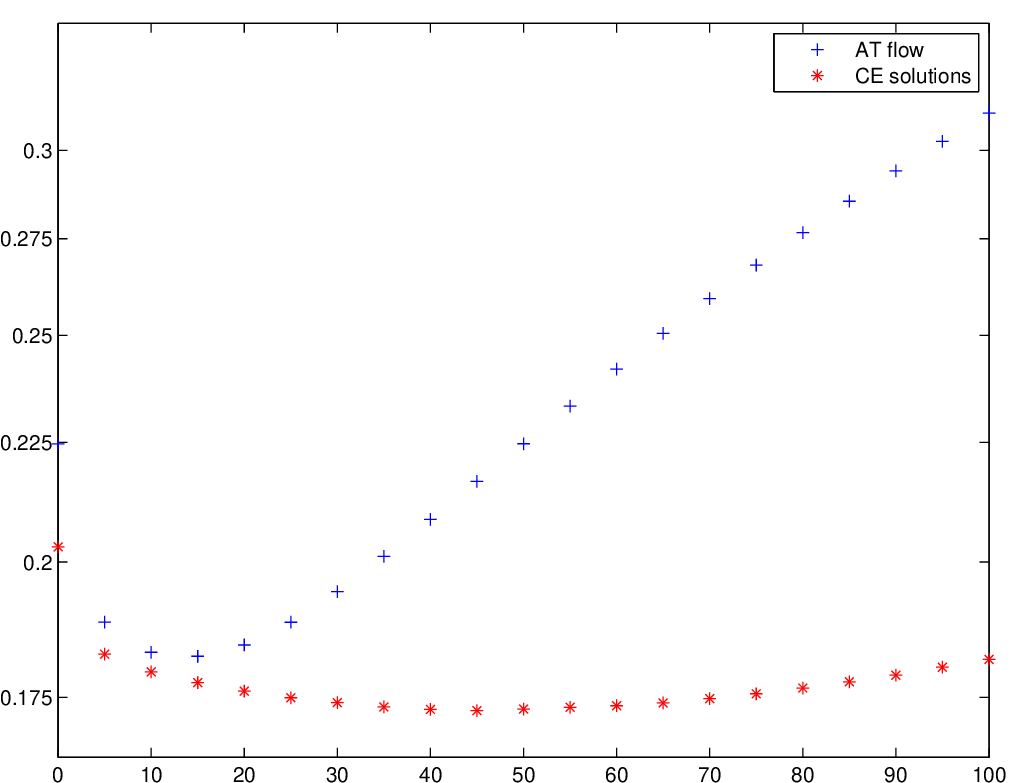}}
\put(200,10){\includegraphics[width=180 pt]{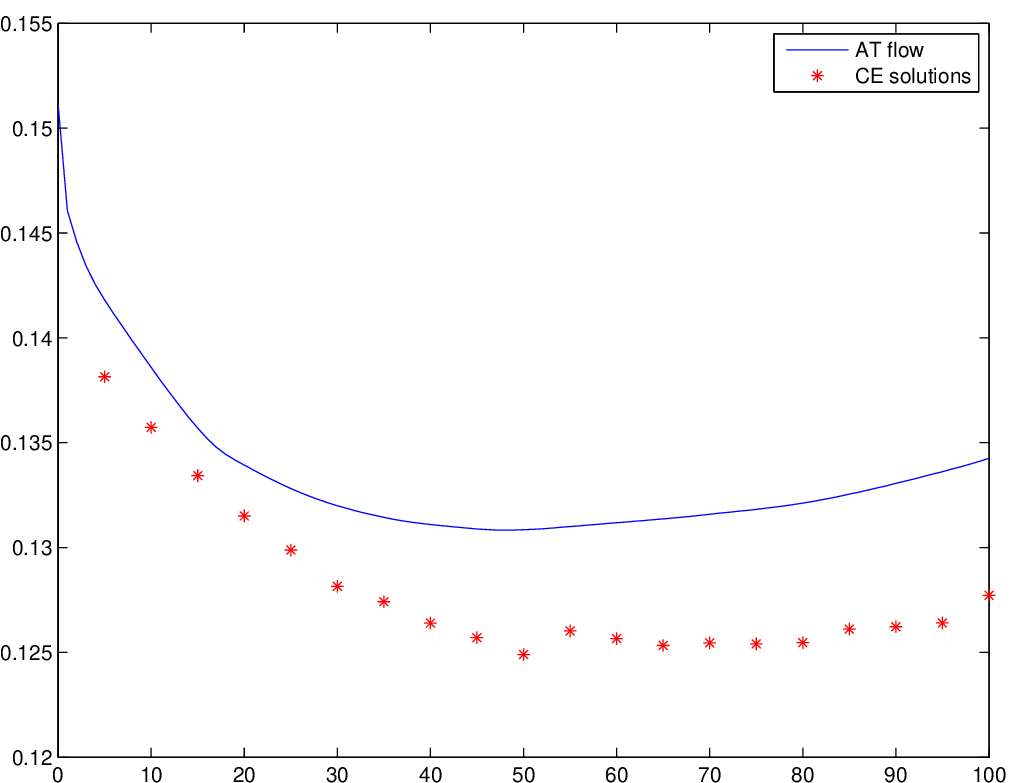}}
\put(150,0){{\scriptsize Number of iterations}}
\put(35,150){{\scriptsize rel. $l^2$-error of CGO sinogram}}
\put(215,150){{\scriptsize rel. $l^1$-error of reconstructed conductivity}}
\end{picture}
\caption{\label{fig:HnL_convergence_plots}Convergence plot of the AT minimized solutions CGO sinograms in relative $l^2$-error to the true measurement (Left) and the relative $l^1$-error of the reconstructions to the true Heart-and-Lungs conductivity (Right).}
\end{figure} 

\clearpage

\subsection{An Industrial Pipe Phantom}\label{sec:pipe}
The second test case is an example from the oil industry. It roughly models a pipe with oil (top layer, conductivity 1.2), water (middle layer, conductivity 2.0), and sand (bottom layer, conductivity 0.3) similar to the test problem of \cite{Astala2011}.

The initial D-bar reconstruction is computed on the same mesh with a truncation radius of $R=6$ and additional noise of $0.01\%$ in the measured DN map. The parameters for the AT flow and contrast enhancement were chosen the same as in the previous example. The {\sc CGO} sinograms were computed with a small radius of $|k|=0.5$ far within the reliable region. 

The minimization of the CGO sinogram error with the proposed algorithm did not converge and hence the algorithm stopped at the maximum number of iterations chosen as $J_{max} = 200$. 

The error in the CGO sinogram has been halved from $17.53\%$ in the original D-bar reconstruction, to $8.73\%$ in the improved reconstruction obtained with the new algorithm (left plot of Figure \ref{fig:Pipe_convergence_plots}). The right plot shows the errors for contrast enhanced versions of the conductivities after the AT flow. Their error increased, after an initial decrease, from $15.09\%$ of the D-bar reconstruction to $15.32\%$ of the reconstruction obtained by the AT flow after 200 iterations and $18.35\%$ of the corresponding contrast enhanced solution.  However, we point out that in practice such a comparison (right plot of Figure~\ref{fig:Pipe_convergence_plots}) is not possible as the true conductivity is unknown, and only the left plot is possible, where in fact we see a clear decrease in the CGO error.


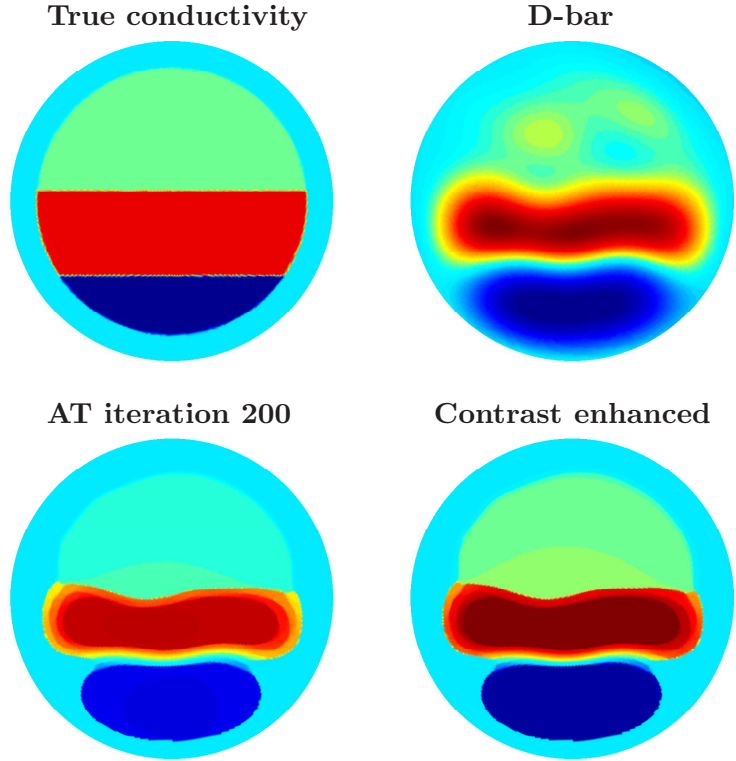
\begin{figure}[ht!]
\centering
\begin{picture}(300,290)
\put(0,150){\includegraphics[width=125 pt]{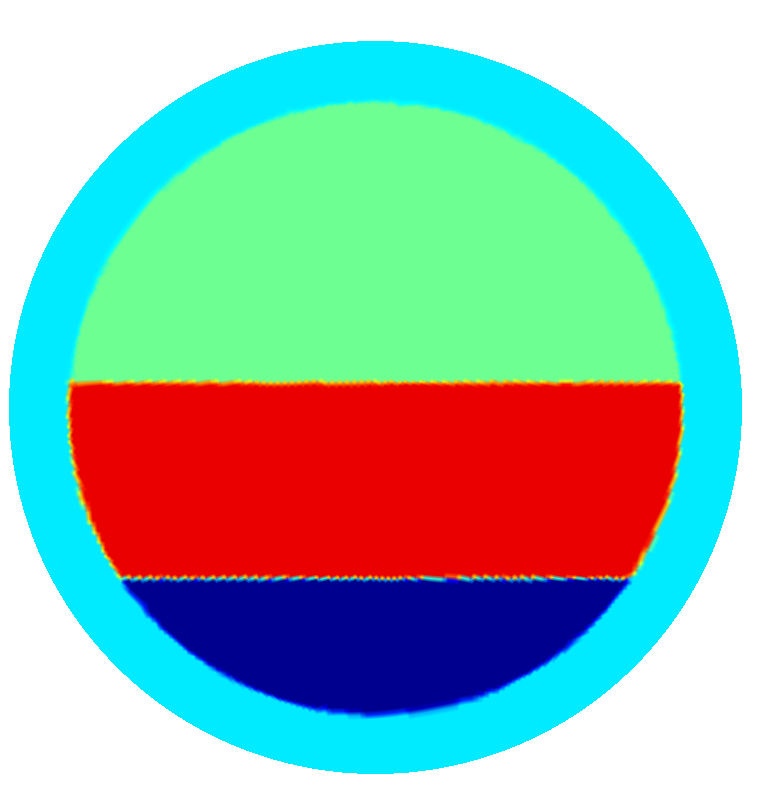}}
\put(150,150){\includegraphics[width=125 pt]{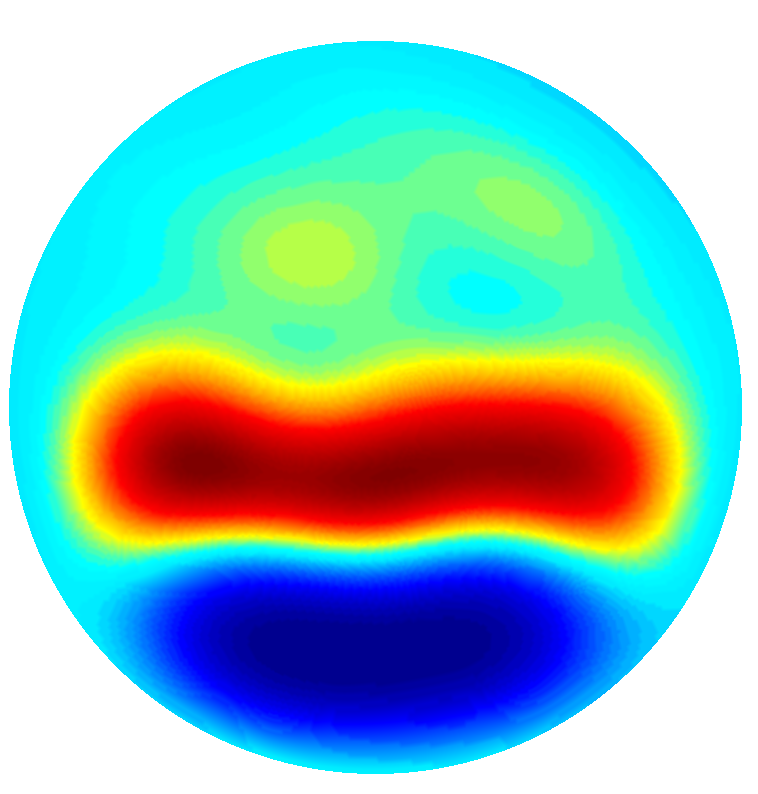}}
\put(15,280){\textbf{True conductivity}}
\put(195,280){\textbf{D-bar}}
\put(0,0){\includegraphics[width=125 pt]{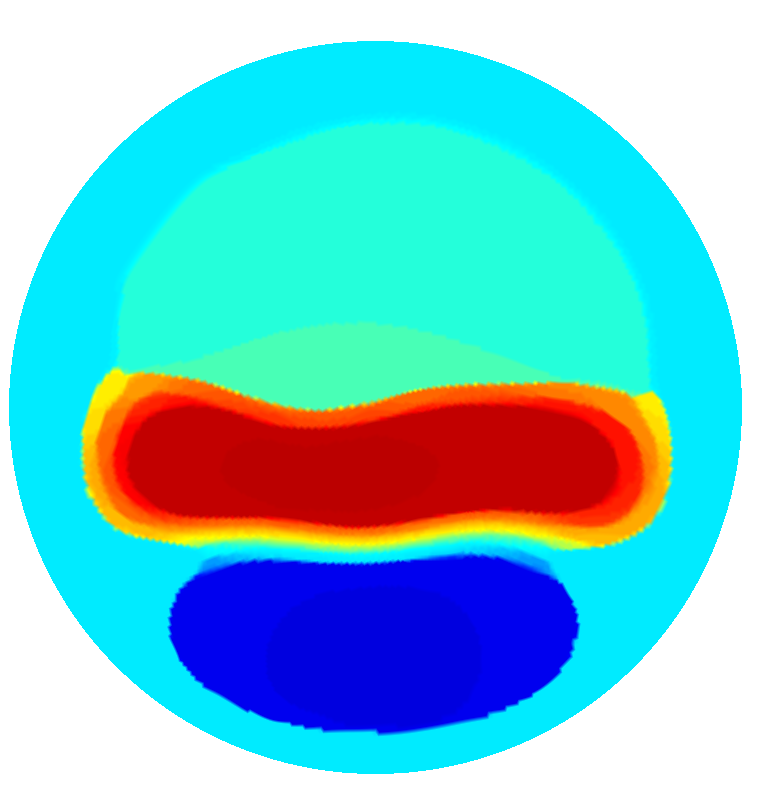}}
\put(150,0){\includegraphics[width=125 pt]{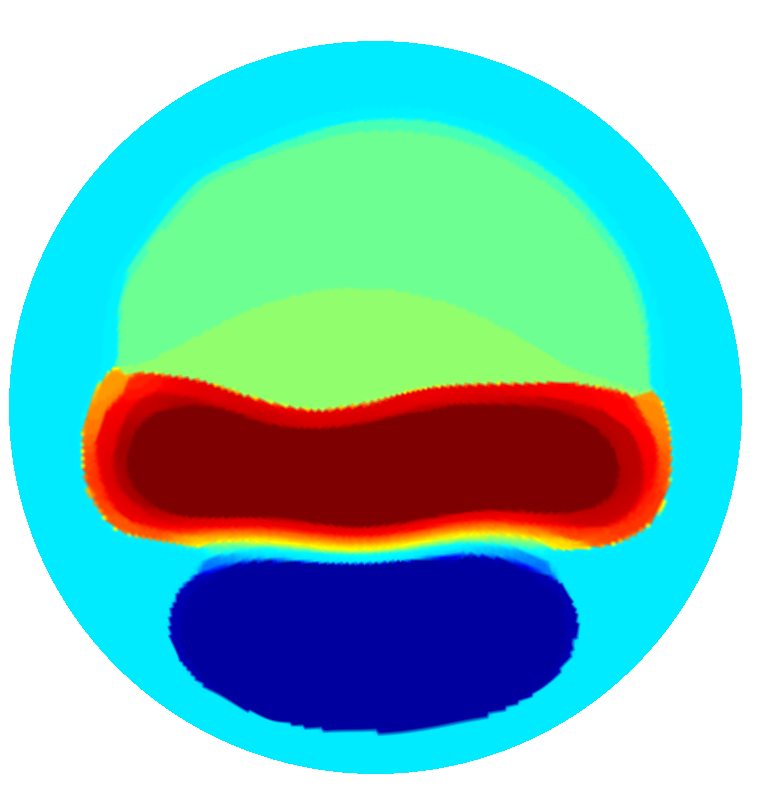}}
\put(15,130){\textbf{AT iteration 200}}
\put(160,130){\textbf{Contrast enhanced}}
\end{picture}
\caption{\label{fig:Pipe_noisy_sum}Illustration of original pipe phantom with the initial D-bar reconstruction (Truncation radius 6) in the top row. Reconstruction after 200 iterations of the AT flow (Bottom left) and the corresponding contrast enhanced solution (Bottom right).}
\end{figure} 
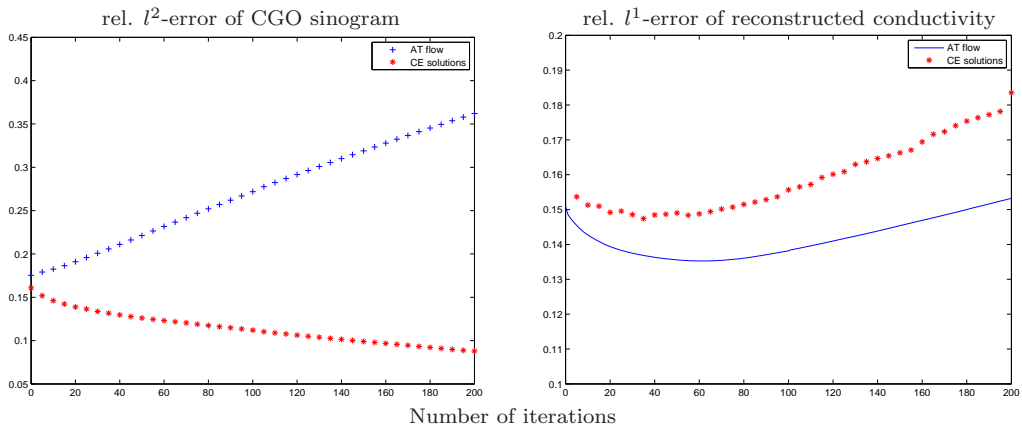
\begin{figure}[ht!]
\centering
\begin{picture}(400,175)
\put(0,10){\includegraphics[width=179 pt]{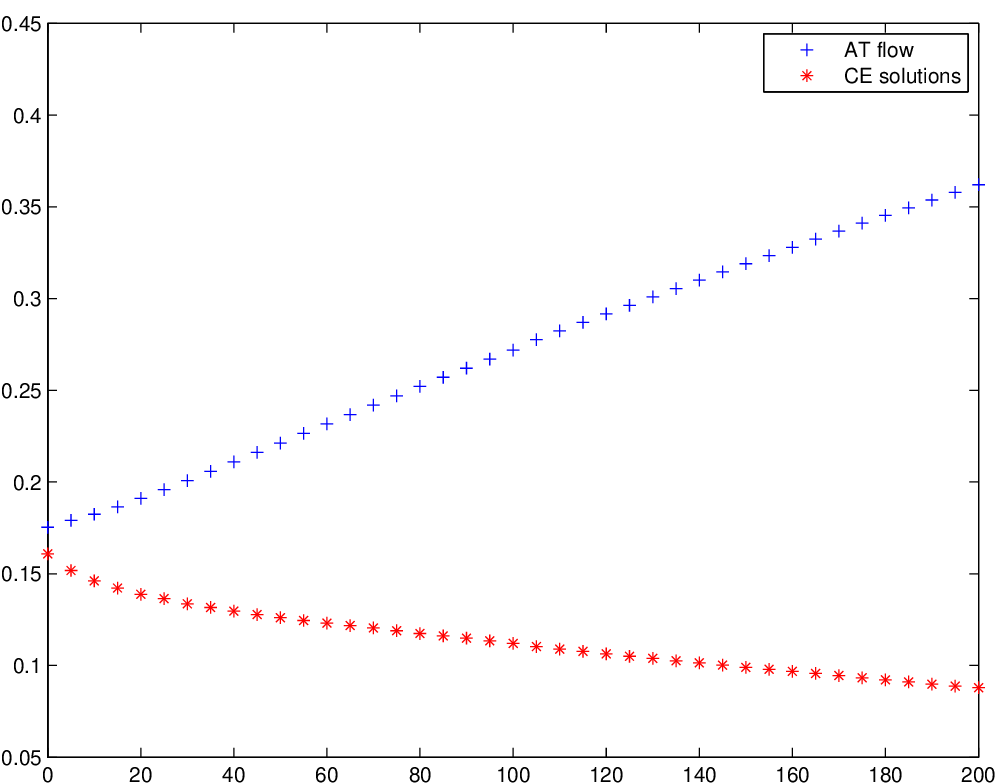}}
\put(200,10){\includegraphics[width=180 pt]{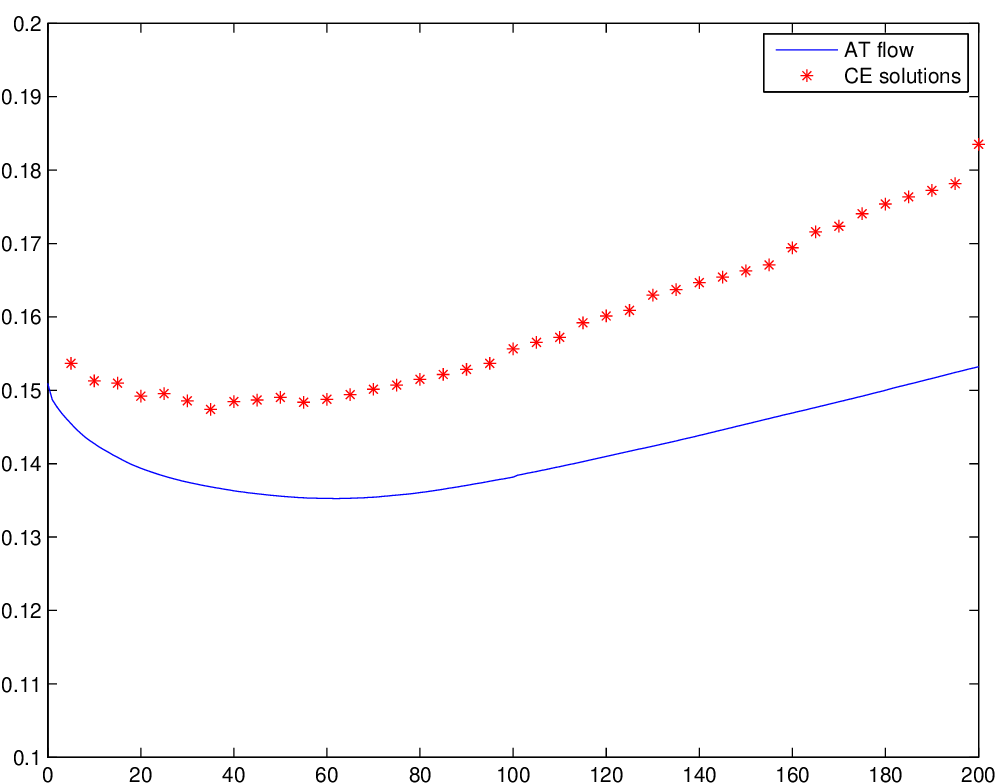}}
\put(150,0){{\scriptsize Number of iterations}}
\put(35,150){{\scriptsize rel. $l^2$-error of CGO sinogram}}
\put(215,150){{\scriptsize rel. $l^1$-error of reconstructed conductivity}}
\end{picture}
\caption{\label{fig:Pipe_convergence_plots} Convergence plot of the AT minimized solutions CGO sinograms in relative $l^2$-error to the true measurement (Left) and the relative $l^1$-error of the reconstructed pipe conductivity (Right).}
\end{figure}

\afterpage{\clearpage}

\section{Discussion}\label{sec:discus}
The heart-and-lungs phantom of Section \ref{sec:heartlungs} is of particular interest, since it represents, an admittedly simplified version of, a major application of electrical impedance tomography in the medical field: monitoring the blood and air flow in a patient's heart and lungs.  With the proposed algorithm we were able to clearly distinguish the left and right lung as well as introduce clear edges, see Figure~\ref{fig:HnL_noisy_sum}.  As discussed above, a level of $0.5\%$ noise is reasonably high for EIT measurements and hence gives a good impression of how the algorithm behaves with noise corrupted data. The original D-bar recovered conductivity (see the middle image in Figure~\ref{fig:HnL_noisy_sum}) has good contrast but lacks sharpness, as is typical for D-bar reconstructions. 	

A summary of the AT flow is illustrated in Figure \ref{fig:HnL_AT_illu}.  Notice that the AT flow gradually divides the two lungs, and that the separated areas converge to constant values.   As seen in the right plot of Figure~\ref{fig:HnL_convergence_plots}, the AT flow (blue line) clearly minimizes the $l^1$-error of the evolved conductivity to the true conductivity.  Furthermore, the behavior of its corresponding {\sc CGO} sinogram error (blue plus sign markers in the left plot), initially decreasing but followed by a sharp increase, reinforces the need for a contrast enhancement step.  Note that the CE solutions have a similar behavior in both {\sc CGO} sinogram error and the reconstruction error of the conductivity, and that the CE evolved conductivities have smaller relative errors to the true conductivity than their non CE counterparts.
\newline

The second test case (Section \ref{sec:pipe}) is an example from the oil industry.  When measuring a pipe with oil (top layer), water (middle layer), and sand (bottom layer) one wants to know how much oil is transported in the pipe, making it important to distinguish the levels of each substance clearly, i.e. their edges.  The new reconstruction has clear and sharp edges dividing the different substances and one can tell easily how much of each substance is present in the pipe. As one can see in Figure \ref{fig:Pipe_noisy_sum}, the structure of the pipe has been reintroduced to the reconstruction, delivering a realistic view of the imaged target.

The right plot in Figure \ref{fig:Pipe_convergence_plots} shows that the error in reconstructed conductivities did not decrease as nicely as in the Heart-and-Lungs phantom. In fact, the \newtext{relative $l^1$ conductivity} reconstruction error increased from $15.09\%$ of the D-bar reconstruction to $18.35\%$ of the contrast enhanced solution of the last iteration 200 in the AT flow. A reason for the higher error can be seen in Figure \ref{fig:Pipe_noisy_sum}, the contrast enhancement increases the conductivity of the top layer as well as the middle layer, which produces a higher error in the middle layer. This suggests that more sensitive models for the contrast enhancement may be needed. Nevertheless, the information contained in the contrast enhanced solution is far more useful for evaluation of the target  and suggests that the CGO sinogram contains more information about the reconstruction's geometry. 

\section{Conclusions}\label{sec:conclude}
\noindent
A novel edge-preserving D-bar method with a data-driven contrast enhancement was introduced and tested on simulated EIT measurement data.  The algorithm works as advertised by both sharpening and enhancing the contrast of the reconstruction, even in the presence of additional noise added to the Dirichlet-to-Neumann boundary measurements.

Key to the approach is the invention of the {\sc CGO} sinogram, a more reliable and geometrically transparent quantity than the DN map.  
\newtext{The {\sc CGO} sinogram provides an important breakthrough towards new uses of this nonlinear data, and based on our findings we are excited what further theoretical analysis on the stability of the CGO sinogram will reveal.}


Our results also have implications outside the realm of D-bar methods. The traditional regularization approach for EIT reconstructions is to find the minimizer of a functional of the form
\begin{equation}\label{tikhonov}
  \|\Lambda_\sigma^\delta-\Lambda_{\sigma^\prime}\|_Y+\alpha \|\sigma^\prime\|_X,
\end{equation}
where the $X$ norm corresponds, for instance, to Tikhonov or Total Variation regularization, and $0<\alpha<\infty$ is a regularization parameter. Replacing the traditional data fidelity term in (\ref{tikhonov}) by an analogous term based on the CGO sinogram leads to the functional
\begin{equation}\label{def:newfunctional}
\|\mathcal{S}_\sigma(\theta,\varphi,r)-\mathcal{S}_{\sigma^\prime}(\theta,\varphi,r)\|_{L^2(\mathbb{T}^2)}^2+\alpha \|\sigma^\prime\|_X,
\end{equation}
where $\mathbb{T}^2$ denotes the two-dimensional torus. Based on the evidence seen in Figures \ref{fig:CGOsinogramframes} and \ref{fig:CGOexample2}, we strongly suspect that using (\ref{def:newfunctional}) would lead to superior reconstructions compared to (\ref{tikhonov}). A similar comment applies to the likelihood distributions used in Bayesian inversion approaches for EIT.

\section*{Acknowledgments}
\noindent
The study was supported by the SalWe Research Program for Mind and Body (Tekes - the Finnish Funding Agency for Technology and Innovation grant 1104/10) and by the Academy of Finland (Finnish Centre of Excellence in Inverse Problems Research 2012--2017, decision number 250215).  

\newpage
\bibliographystyle{siam}
\bibliography{Inverse_problems_references_HHS}

\begin{thebibliography}{10}

\bibitem{Alicandro1998}
{\sc R.~Alicandro, A.~Braides, and J.~Shah}, {\em Approximation of non-convex
  functionals in {GBV}}, 1998.

\bibitem{Ambrosio1990}
{\sc L.~Ambrosio and V.~M. Tortorelli}, {\em Approximation of functionals
  depending on jumps by elliptic functionals via {$\Gamma$}-convergence},
  Communications on Pure and Applied Mathematics, 43 (1990), pp.~999--1036.

\bibitem{Astala2011}
{\sc K.~Astala, J.~Mueller, L.~P{\"a}iv{\"a}rinta, A.~Per{\"a}m{\"a}ki, and
  S.~Siltanen}, {\em Direct electrical impedance tomography for nonsmooth
  conductivities}, Inverse Problems and Imaging, 5 (2011), pp.~531--549.

\bibitem{Astala2006}
{\sc K.~Astala and L.~P\"aiv\"arinta}, {\em A boundary integral equation for
  {C}alder\'on's inverse conductivity problem}, in Proc. 7th Internat.
  Conference on Harmonic Analysis, Collectanea Mathematica, 2006.

\bibitem{Astala2006a}
{\sc K.~Astala and L.~P{\"a}iv{\"a}rinta}, {\em {C}alder\'on's inverse
  conductivity problem in the plane}, Annals of Mathematics, 163 (2006),
  pp.~265--299.

\bibitem{Bikowski2010a}
{\sc J.~Bikowski, K.~Knudsen, and J.~L. Mueller}, {\em Direct numerical
  reconstruction of conductivities in three dimensions using scattering
  transforms}, Inverse Problems, 27 (2011), p.~015002.

\bibitem{Brown1996}
{\sc R.~M. Brown}, {\em Global uniqueness in the impedance imaging problem for
  less regular conductivities}, SIAM Journal on Mathematical Analysis, 27
  (1996), pp.~1049--1056.

\bibitem{Calder'on1980}
{\sc A.-P. {C}alder{\'o}n}, {\em On an inverse boundary value problem}, in
  Seminar on {N}umerical {A}nalysis and its {A}pplications to {C}ontinuum
  {P}hysics ({R}io de {J}aneiro, 1980), Soc. Brasil. Mat., Rio de Janeiro,
  1980, pp.~65--73.

\bibitem{Chambolle1995a}
{\sc A.~Chambolle}, {\em Image segmentation by variational methods: Mumford and
  shah functional and the discrete approximations}, SIAM Journal on Applied
  Mathematics, 55 (1995), pp.~827--863.

\bibitem{Cheney1999}
{\sc M.~Cheney, D.~Isaacson, and J.~C. Newell}, {\em Electrical impedance
  tomography}, SIAM Review, 41 (1999), pp.~85--101.

\bibitem{Cornean2006}
{\sc H.~Cornean, K.~Knudsen, and S.~Siltanen}, {\em Towards a {$d$}-bar
  reconstruction method for three-dimensional {\sc {eit} }}, Journal of Inverse
  and Ill-Posed Problems, 14 (2006), pp.~111--134.

\bibitem{INCR}
{\sc I.~N.~R. Council}, {\em Dielectric properties of body tissues}, 2013.

\bibitem{DeGiorgi1989}
{\sc E.~De~Giorgi, M.~Carriero, and A.~Leaci}, {\em Existence theorem for a
  minimum problem with free discontinuity set}, Archive for Rational Mechanics
  and Analysis, 108 (1989), pp.~195--218.

\bibitem{Delbary2011}
{\sc F.~Delbary, P.~C. Hansen, and K.~Knudsen}, {\em Electrical impedance
  tomography: 3{D} reconstructions using scattering transforms}, Applicable
  Analysis, 0 (0), pp.~1--19.

\bibitem{Erdem2009}
{\sc E.~Erdem and S.~Tari}, {\em Mumford-shah regularizer with contextual
  feedback}, Journal of Mathematical Imaging and Vision, 33 (2009), pp.~67--84.

\bibitem{Faddeev1966}
{\sc L.~D. Faddeev}, {\em Increasing solutions of the {S}chr{\"o}dinger
  equation}, Soviet Physics Doklady, 10 (1966), pp.~1033--1035.

\bibitem{Finkel2003}
{\sc D.~E. Finkel}, {\em Direct optimization algorithm user guide}, tech. rep.,
  Center for Research in Scientific Computation, North Carolina State
  University, 2003.

\bibitem{Francini2000}
{\sc E.~Francini}, {\em Recovering a complex coefficient in a planar domain
  from {D}irichlet-to-{N}eumann map}, Inverse Problems, 16 (2000),
  pp.~107--119.

\bibitem{Hamilton2012}
{\sc S.~Hamilton, C.~Herrera, J.~L. Mueller, and A.~VonHerrmann}, {\em A direct
  {D-bar} reconstruction algorithm for recovering a complex conductivity in
  {2-D}}, Inverse Problems, 28 (2012), p.~095005.

\bibitem{HM12_NonCirc}
{\sc S.~J. Hamilton and J.~L. Mueller}, {\em Direct {EIT} reconstructions of
  complex admittivities on a chest-shaped domain in {2-D}}, IEEE Transactions
  on Medical Imaging, 32 (2013), pp.~757--769.

\bibitem{Harhanen2014}
{\sc L.~{Harhanen}, N.~{Hyv{\"o}nen}, H.~{Majander}, and S.~{Staboulis}}, {\em
  {Edge-enhancing reconstruction algorithm for three-dimensional electrical
  impedance tomography}}, ArXiv e-prints, arXiv:1406.1279 (2014).

\bibitem{Helin2011}
{\sc T.~Helin and M.~Lassas}, {\em Hierarchical models in statistical inverse
  problems and the mumford--shah functional}, Inverse problems, 27 (2011),
  p.~015008.

\bibitem{Isaacson2006}
{\sc D.~Isaacson, J.~Mueller, J.~Newell, and S.~Siltanen}, {\em Imaging cardiac
  activity by the {D}-bar method for electrical impedance tomography},
  Physiological Measurement, 27 (2006), pp.~S43--S50.

\bibitem{Jones1993}
{\sc D.~R. Jones, C.~D. Perttunen, and B.~E. Stuckman}, {\em Lipschitzian
  optimization without the lipschitz constant}, Journal of Optimization Theory
  and Applications, 79 (1993), pp.~157--181.

\bibitem{Jung2011}
{\sc M.~Jung, X.~Bresson, T.~F. Chan, and L.~A. Vese}, {\em Nonlocal
  {M}umford-{S}hah regularizers for color image restoration}, IEEE Trans. Image
  Process., 20 (2011), pp.~1583--1598.

\bibitem{Knudsen2003}
{\sc K.~Knudsen}, {\em A new direct method for reconstructing isotropic
  conductivities in the plane}, Physiological Measurement, 24 (2003),
  pp.~391--403.

\bibitem{Knudsen2009}
{\sc K.~Knudsen, M.~Lassas, J.~Mueller, and S.~Siltanen}, {\em Regularized
  {D}-bar method for the inverse conductivity problem}, Inverse Problems and
  Imaging, 3 (2009), pp.~599--624.

\bibitem{Knudsen2004a}
{\sc K.~Knudsen and A.~Tamasan}, {\em Reconstruction of less regular
  conductivities in the plane}, Communications in Partial Differential
  Equations, 29 (2004), pp.~361--381.

\bibitem{Mueller2003}
{\sc J.~Mueller and S.~Siltanen}, {\em Direct reconstructions of conductivities
  from boundary measurements}, SIAM Journal on Scientific Computing, 24 (2003),
  pp.~1232--1266.

\bibitem{Mueller2012}
{\sc J.~Mueller and S.~Siltanen}, {\em Linear and Nonlinear Inverse Problems
  with Practical Applications}, vol.~10 of Computational Science and
  Engineering, SIAM, 2012.

\bibitem{Mumford1985}
{\sc D.~Mumford and J.~Shah}, {\em Boundary detection by minimizing
  functionals}, in IEEE Conference on Computer Vision and Pattern Recognition,
  1985.

\bibitem{Music2013}
{\sc M.~Music, P.~Perry, and S.~Siltanen}, {\em Exceptional circles of radial
  potentials}, Inverse Problems, 29 (2013), p.~045004.

\bibitem{Nachman1996}
{\sc A.~I. Nachman}, {\em Global uniqueness for a two-dimensional inverse
  boundary value problem}, Annals of Mathematics, 143 (1996), pp.~71--96.

\bibitem{Perona1990}
{\sc P.~Perona and J.~Malik}, {\em Scale-space and edge detection using
  anisotropic diffusion}, IEEE Transactions on Pattern Analysis and Machine
  Intelligence, 12 (1990), pp.~629--639.

\bibitem{Ramlau2007}
{\sc R.~Ramlau and W.~Ring}, {\em A mumford--shah level-set approach for the
  inversion and segmentation of x-ray tomography data}, Journal of
  Computational Physics, 221 (2007), pp.~539--557.

\bibitem{Rondi2001}
{\sc L.~Rondi and F.~Santosa}, {\em Enhanced electrical impedance tomography
  via the {M}umford-{S}hah functional}, ESAIM Control Optim. Calc. Var., 6
  (2001), pp.~517--538.

\bibitem{Shah1996}
{\sc J.~Shah}, {\em A common framework for curve evolution, segmentation and
  anisotropic diffusion}, in IEEE Conference on Computer Vision and Pattern
  Recognition, IEEE, 1996, pp.~136--142.

\bibitem{Siltanen2000}
{\sc S.~Siltanen, J.~Mueller, and D.~Isaacson}, {\em An implementation of the
  reconstruction algorithm of a. nachman for the 2-d inverse conductivity
  problem}, Inverse Problems, 16 (2000), pp.~681--699.

\bibitem{Sylvester1987}
{\sc J.~Sylvester and G.~Uhlmann}, {\em A global uniqueness theorem for an
  inverse boundary value problem}, Annals of Mathematics, 125 (1987),
  pp.~153--169.

\bibitem{Weickert1998}
{\sc J.~Weickert}, {\em Anisotropic Diffusion in Image Processing}, Teubner
  Stuttgart, 1998.

\end{thebibliography}

\end{document}